\documentclass[12pt]{article}
\usepackage{amsmath,amsfonts,latexsym}

\usepackage{amsfonts}

\def\C{{\Bbb C}}

\def\Z{{\Bbb Z}}

\def\N{{\Bbb N}}

\begin{document}
\begin{center}
{\LARGE{\bf Classification of Irreducible integrable modules for
toroidal Lie-algebras with finite dimensional weight spaces}}  \\ [5mm]
{\bf S. Eswara Rao}
\end{center}
\vskip 5mm
\centerline{Dedicated to Professor P. Somaraju with admiration and gratitude.}
\vskip 1cm

The study of Maps $(X,G)$, the group of polynomial maps of a complex
algebraic variety $X$ into a complex algebraic group $G$, and its
representations is only well developed in the case that $X$ is a
complex torus $\C^*$.  In this case Maps $(X, G)$ is a loop group and
the corresponding Lie-algebra Maps $(X, \stackrel{\circ}{\cal G})$ is
the loop algebra $\C [t, t^{-1}] \otimes \stackrel{\circ}{\cal G}$.
Here the representation comes to life only after one replaces Maps
$(X, \stackrel{\circ}{\cal G})$ by its universal central extension,
the corresponding affine Lie-algebra.  One then obtains the well
known theory of highest weight  modules, vertex representations,
modular forms and character theory and so on.

The next easiest case is presumably the case of $n$ dimensional torus
$(\C^*)^n$.  So we consider the universal central
extension of $\stackrel{\circ}{\cal G}\otimes \C
[t_{1}^{\pm}, \cdots, t_{{n}}^{\pm}]$ which is
referred to as the toroidal Lie-algebra $\tau$ in $[EM]$ and $[MEY]$.  The
most interesting modules are the integrable modules (where the real
root space acts locally nilpotently (see section 2)), as they lift to
the corresponding group.  Unlike the affine case where the central
extension is one dimensional, the toroidal case has infinite
dimensional centre which makes the theory more complicated.  For the first
time a large number of integrable (reducible) modules for toroidal
Lie-algebras (simply laced case) have been constructed through the
use of vertex operators in [EM] and [MEY].

In this paper we construct two classes of (Examples (4.1) and (4.2))
of irreducible integrable modules for toroidal Lie-algebras with
finite dimensional weight spaces.  In sections 4 and 5 we prove that
any irreducible integrable module where part of the centre acts
non-trivially are the ones given in Example (4.2) upto an
automorphism of $\tau$.  We have proved in [E3] that the only
modules with the above property where center acts trivially are the
one given in Example 4.1 (see Remark 5.5).  In this case a similar
classification is obtained in [YY] with a stronger assumption on the
weight spaces. These results in the case $n=1$ are due to [C] and [CP].

In section 1, after recalling the construction of non-twisted affine
Lie-algebras from [K], we establish the necessary terminology for
root systems, non-degenerate billinear form and the Weyl group for
toroidal Lie-algebras.  In section 2 we recall the actual definition of
toroidal Lie-algebras and integrable modules and prove some standard
facts about the Weyl group and weight systems (Lemma 2.3). Then we
prove that an irreducible integrable module for $\tau$ with finite
dimensional weight spaces has a highest weight vector in the following
sense.  Let ${\cal G}_{af}= \stackrel{\circ}{\cal G}\otimes \C
[t_1, t_{1}^{-1}] \oplus \C C_1 \oplus \C d_1$ be an affine
Lie-algebra and let ${\cal G}_{af}= N^- \oplus h^1
\oplus N^+$ be the standard decomposition.  Then we prove that if
some zero degree central operator acts non-trivially there exists
a vector killed by $N^+ \otimes A_{n-1}$ (Proposition (2.4) and
(4.8) after twisting the module upto an
automorphism).  Let $\stackrel{\circ}{\cal G}=n^- \oplus h^\circ
\oplus n^+$ and suppose all zero degree central operators act
trivially then we prove that there is a vector killed by $n^+ \otimes
A_n$ (Proposition 2.12).  These two results may be seen as
generalization of Theorem (2.4) of [C].

In the section 3 we define graded and non-graded highest weight
modules in the generality of loop Kac-Moody Lie-algebras and prove
(Proposition 3.5) that there is a one-one correspondence between
the graded and non-graded cases.  The problem now reduces to the
classification of irreducible integrable highest weight  modules
(non-graded) (Lemma 3.6).  We prove (Remark 3.9) that any such
module is actually a module for ${\cal G}_{af} \otimes A_{n-1} /
I$ for a co-finite ideal $I$.  (we are in the case where some zero
degree center acts non-trivially).  In a combinatorial Lemma
(3.11) which is of independent interest we prove that such a
Lie-algebra is isomorphic to $\oplus {\cal G}_{af}$ (direct sum of
finitely many copies of ${\cal G}_{af})$.  Then it is very
standard to classify irreducible integrable highest weight modules
for $\oplus {\cal G}_{af}$.

If the zero degree center acts trivially then the full center
should act trivially (Proposition 4.13)  where we use an
interesting result (Proposition 4.12) on Hisenberg Lie-algebras
due to Futorny [F].  In this case the classification is given in
[E3] (see Remark 5.5).

We prove in Lemma (4.6) in the generality of irreducible modules for
$\tau$ with finite dimensional weight spaces that most of the center
acts trivially.  In fact in each graded component of the center at
most one dimensional space acts non-trivially.

In other papers [BB] and [E4] a more general toroidal Lie-algebra is
considered by adding an infinite set of derivations.  They have
constructed integrable irreducible (highest weight) modules for
toroidal Lie-algebra where almost all of the ``center'' acts
non-trivially.
It remains to be seen which toroidal Lie-algebra admits an interesting
representation theory in general.  More general Lie-algebras called
extended affine Lie-algebras (EALA), of which toroidal Lie-algebras
are prime examples, are studied extensively.  See for example [BGK]
and the references therein.
\section*{1 Section}

We first recall the some notation of the theory of affine Lie-algebras from
chapter 6 of [K].  We always denote $\N$ the non-negative
integers, $\Z$ the integers and $\Z_+$ the positive integers.  All our
vector spaces are over the complex numbers $\C$.

Let $\stackrel{\circ}{\cal G}$ be a simple finite dimensional
Lie-algebra.  Let $\stackrel{\circ}{h} $ be a Cartan subalgebra of
finite dimension $d$.  Let $A= (a_{ij})$ be the Cartan matrix of the
corresponding non-twisted affine Lie-algebra.  Let $a_0, a_1, \cdots$
$a_d$ be the numerical labels of the Dynkin diagram $S (A)$.  Let
$a_0^\vee, \cdots, a_d^\vee$ be the numerical labels of the Dynkin diagram
$S(A^T)$.  Let $b_i^{-1} = a_i^\vee$.  We know from [K] that for a
non-twisted affine Lie-algebra that $a_0=b_0=1$.  Write $A=Diag (a_0
b_0, a_1 b_1, \cdots, a_d b_d)B$ where $B$ is a symmetric matrix.  Let
$h$ be the Cartan subalgebra spanned by $\alpha_0^\vee, \alpha_1^\vee,
\cdots, \alpha_d^\vee, d_0$ such that $\stackrel{\circ}{h}$ is spanned
by $\alpha_1^\vee, \cdots, \alpha_d^\vee$.  Define $\alpha_i \in h^*$ by
$\alpha_i (\alpha_j^\vee) = a_{ji}, \alpha_i (d_0) = 0$ for $1 \leq i
\leq d$ and $\alpha_0 (d_0)=1$.  Define $w$ in $h^*$ such that
$w (\alpha_i^\vee) = \delta_{i0}$ and $w (d_0)=0$.  Then
$h^*$ is spanned by $\alpha_0, \alpha_1, \cdots \alpha_d,w$
for dimension reasons.  Define nondegenerate symmetric
billinear form on $h$ by $(\alpha_i^\vee, \alpha_j^\vee) = a_j b_j a_{ij},
(\alpha_i^\vee, d_0) = \delta_{i0}$ and $(d_0, d_0) =0$.  Similarly
$h^*$ carries a billinear form such that $(\alpha_i, \alpha_j) =
b_i^{-1} a_i^{-1} a_{ij}, (\alpha_i, w) = \delta_{i0}$ and
$(w,w) =0$.  As in [K] we normalize the form on $h^*$
by $(\beta, \beta) =2$ where $\beta$ is the maximal root of
$\stackrel{\circ}{{\cal G}}$. Define
$$
\begin{array}{lll}
C&=& \sum_{i=1}^{d} b_i^{-1} \alpha_i^\vee + \alpha_0^\vee \\
\delta &=& \sum_{i=1}^{d} a_i \alpha_i + \alpha_0
\end{array}
$$
Then $\delta- \alpha_0 = \beta$.  Let $\beta^\vee = C - \alpha_0^\vee$.
>From [K] we have the following.

$(\delta, \alpha_i) = 0, (\delta, \delta) =0, (\delta, w)=1,
(C, \alpha_i^\vee)=0, (C,C)=0, (C, d_0)=1, \delta (\alpha_i^\vee)=0,
\alpha_i (C) =0$ and $w (C)=1$.

\paragraph*{(1.1) Lemma}  For $\lambda \in h^*$
$$
\begin{array}{lll}
(1) ~~~ (\lambda, \alpha_i) &=& b_i^{-1} a_i^{-1} \lambda
(\alpha_i^\vee), 0 \leq i \leq d \\
(2) ~~~ (\lambda, \beta) &=& \lambda (\beta^\vee)
\end{array}
$$
\paragraph*{Proof:}  It is easy to see.

Define the $i$th fundamental reflection $r_i$ on $h^*$ such that $r_i
\lambda = \lambda - \lambda (\alpha_i^\vee) \alpha_i, 0 \leq i \leq d$.
Let $W_{af}$ be the Weyl group
generated by $r_i$.  Then  (,) on $h^*$ is $W_{af}-$invariant.

We need the following simple lemma.  Let $\stackrel{0}{\Lambda}$ be the root
lattice of $\stackrel{\circ}{{\cal G}}$ spanned by $\alpha_1,
\alpha_2, \cdots, \alpha_d$.  For $\lambda, \mu \in
\stackrel{\circ}{h}^*$ define $\lambda  \frac{<}{\circ} \mu$
if $\mu - \lambda =
\displaystyle{\sum_{i=1}^{d}} n_i \alpha_i$ for $n_i \in \N$.  Let
$\stackrel{0}{\Lambda}^+$ be the set of weights in
$\stackrel{\circ}{h}^*$ such
that $\lambda (\alpha_i^\vee) \in \N$ for $1 \leq i \leq d$.  A weight
$\lambda \in \stackrel{0}{\Lambda}^+$ is called miniscule if $\mu
\frac{<}{\circ}
\lambda, \mu \in \stackrel{0}{\Lambda}^+$ implies $\mu= \lambda$.

\paragraph*{(1.2) Lemma}  Let $\lambda \in \stackrel{\circ}{h}^*$ be
miniscule then $\lambda (\beta^\vee) = 0 $ or $1$.

\paragraph*{Proof} Follows from
exercise 13 of Chapter III of [H].  Just note that $\lambda$ is
dominant integral.

We will now generalize this for toroidal algebras ([EM], [MEY]).
Fix a positive integer $n$.  Let $\underline{h}$ be the $2n+d$
dimensional vector space spanned by $\alpha_1^\vee, \alpha_2^\vee,
\cdots \alpha^\vee_{d+n}, d_1, \cdots, d_n$.  Let $\tilde{A} =
(A_{ij})$ be a matrix of order $n+d$ such that removal of $n-1$
rows and the corresponding columns in the last $n$ rows and $n$
columns should give $A$. We will describe the matrix explicitly.
Recall that $A=(a_{i,j})_{0\le i,j\le d}$ is the affine matrix.
Then $A_{d+i,j}=a_{0,j},1\le i\le n, 1\le j\le d$;
$A_{j,d+i}=a_{j,0},1\le i\le n, 1\le j\le d$; $A_{d+i,d+j}=2,1\le
i,j \le n$. Let $\tilde{D} =$ diag $(a_1 b_1, a_2 b_2, \cdots a_d
b_d, 1, \cdots 1)$ be a matrix of order $d+n$. Then $\tilde{A}=
\tilde{D} \tilde{B}$ where $\tilde{B}$ is symmetric matrix.

Define $\alpha_i \in \underline{h}^*$ such that $\alpha_i
(\alpha_j^\vee) = a_{ji}$ for $1 \leq i, j \leq d+n, \alpha_i (d_j)=0$
for $1 \leq i \leq d, 1 \leq j \leq n, \alpha_{d+i} (d_j)=\delta_{ij}$
for $ 1 \leq i, j \leq n$.  Define $w_i \in \underline{h}^* (1
\leq i \leq n)$ by $w_i (\alpha_j^\vee) =0$ for $1 \leq j \leq d$
and $w_i (\alpha^\vee_{j+d}) = \delta_{ij}$ for $1 \leq j \leq n, w_j
(d_i)=0$.
Then we will see that $\alpha_1, \cdots \alpha_{d+n}, w_1,
\cdots w_n$ is a  basis of $\underline{h}^*$.  Let $a_{d+i} =
b_{d+j}= 1$ for $1 \leq i, j \leq n$.  Define a symmetric billinear
form $\underline{h}^*$ by $(\alpha_i, \alpha_j) = b_i^{-1}
a_i^{-1} a_{ij}, (\alpha_i, w_{j}) =0 \ 1 \leq i \leq d,
(\alpha_{d+i}, w_j ) = \delta_{ij}$ for $1
\leq i, j \leq n, (w_i, w_j)=0$.

Define for $1 \leq j \leq n$
$$
\begin{array}{lll}
\delta_j &=& \sum_{i=1}^{d} a_i \alpha_i + \alpha_{d+j} = \beta+
\alpha_{j+d} \\
C_j &=& \sum_{i=1}^{d} b_i^{-1} \alpha_i^\vee + \alpha_{d+j}^\vee =
\beta^\vee + \alpha^\vee_{j+d}.
\end{array}
$$

As in the affine case we have the
following $(\alpha_i, \delta_j) =0, (\delta_i, \delta_j)=0,
(\alpha_i^\vee, C_j) =0, (C_i, C_j)=0, \  \delta_i (\alpha_i^\vee)=0, \delta_j
(C_i) =0, \alpha_i (C_j) =0, \delta_j (d_i)= \alpha_{d+j} (d_i)=
\delta_{ij}$ and $w_i (C_j) = \delta_{ij}$.  For $\underline{m}
= (m_1, \cdots m_n) \in \Z^n$ define $\delta_{\underline{m}} = \sum
m_i \delta_i$ and note that $(\delta_{\underline{m}},
\delta_{\underline{m}}) =0$.  $\delta_{\underline{m}}$ are
called null roots.

\paragraph*{(1.3) ~ Root and Co-roots.}  Let
$\stackrel{\circ}{\triangle}$ be the
finite root system.  Let $\triangle = \{\alpha+ \delta \mid \alpha
\in \stackrel{\circ}{ \triangle} \cup \{0\}, \delta$ null root $\}$ \\
$\gamma \in \triangle$ is called real root if $(\gamma,\gamma) \neq 0$.
For $\alpha \in \stackrel{\circ}{\triangle}$ define co-root
$$
\begin{array}{lll}
\alpha^\vee &=& \sum m_i \frac{\mid \alpha_i \mid^2}{\mid \alpha \mid^2}
\alpha_i^\vee \ {\rm for} \ \alpha= \sum_{i=1}^{d} m_i \alpha_i, \alpha
\in \stackrel{\circ}{\triangle}\\
\end{array}
$$
For $\gamma = \alpha + \delta_{\underline{m}}, \alpha \in
\stackrel{\circ}{\triangle} $
Define $\gamma^\vee = \alpha^\vee+ \frac{2}{\mid \alpha \mid^2} \sum m_j
C_j$.

\paragraph*{(1.4) Lemma} For $\lambda \in \underline{h}^*$

$$
\begin{array}{lll}
(1) ~~~ (\lambda, \alpha_i) &=& b_i^{-1} a_i^{-1} \lambda
(\alpha_i^\vee), 0 \leq i \leq d+n \\
(2) ~~~ (\lambda, \alpha) &=& \frac{\mid \alpha \mid^2}{2} \lambda
(\alpha^\vee) , \alpha \epsilon \triangle^{re}
\end{array}
$$
\paragraph*{(1.5) ~ Weyl Group} For a real root $\gamma$ define a
reflection $r_{\gamma}$ on $\underline{h}^*$ by
$$r_\gamma (\lambda) = \lambda - \lambda (\gamma^\vee) \gamma.$$
>From above one can see that it is a  reflection and (,) is
$W-$invariant where $W$ is the Weyl group generated by
$r_\gamma, \gamma$ real.

Now the action of the linear functionals $\alpha_1, \cdots, \alpha_d,
\delta_1, \cdots, \delta_{n}, w_1, w_2, \cdots w_n$ on the
basis $\alpha^\vee_1, \cdots, \alpha^\vee_d, C_1, \cdots C_n, d_1,
\cdots d_n$ is given by the matrix.
$$
\begin{pmatrix}
 \stackrel{\circ}{A} \ 0 \ 0 \cr 0 \ 0 \ I  \cr 0 \ I \ 0
\end{pmatrix}
$$
which is invertible.  Here $\stackrel{\circ}{A}$ is the finite Cartan
matrix.  In particular the linear functions are a basis of
$\underline{h}^*$.

\paragraph*{{\bf 2 ~ Section}} ([EM], [MEY])  Let $\stackrel{\circ}{{\cal
G}}$ be finite dimensional simple Lie-algebra.  Let $n\ge2$ be a
positive integer.  Let $A =A_n = \C [t_1^{\pm}, \cdots t_n^{\pm}]$
be a Laurent polynomial in $n$ commuting variables.  For
$\underline{m}= (m_1, \cdots, m_n) \in \Z^n$ let
$t^{\underline{m}} = t_1^{m_n} \cdots t_n^{m_n}$.  For any vector
space $V$ let $V_A = V \otimes A$ and let $v (\underline{m})= v
\otimes t^{\underline{m}}$ for $v \in V$.  Let ${\cal Z} =
\Omega_A / d_A$ be the module of differentials. That is ${\cal Z}$
is spanned by vectors $t^{\underline{m}} K_i , \leq i \leq n,
\underline{m} \in \Z^n$ by the relation.

\paragraph*{(2.1)} $\sum m_i t^{\underline{m}} K_i =0$

Let $\tau = \stackrel{\circ}{\cal G}_A \oplus {\cal Z} \oplus D$
where $D$ is a linear span of $d_1, \cdots, d_n$.  We will now define
Lie structure on $\tau$ called toroidal Lie-algebra

\paragraph*{(2.2)}
(1) $ [X (\underline{r}), Y (\underline{s})] = [X,Y]
(\underline{r}+\underline{s})+ (X,Y) d (t^{\underline{r}})
t^{\underline{s}}$ where $d (t^{\underline{r}}) t^{\underline{s}} =
\sum r_i t^{\underline{r}+ \underline{s}}K_i$ and (,) is
non-degenerate $\stackrel{\circ}{\cal G}$-invariant symmetric billinear
form on $\stackrel{\circ}{\cal G}$ whose
restriction to $\stackrel{\circ}{h}$ is given in earlier section.

(2)${\cal Z}$ is central in $\stackrel{\circ}{\cal G}_A \oplus {\cal
Z}$

(3) $[d_i, X (\underline{r}) ] = r_i X (\underline{r}) $

(4) $[d_i, t^{\underline{m}} K_j ] = m_i t^{\underline{m}} K_j$.

Let $\underline{h} $ be the span of $\stackrel{\circ}{h}$ and $K_1,
\cdots K_n, d_1, \cdots d_n$.  Clearly $\underline{h}$ can be
identified with $\underline{h}$ defined in section 1 with $K_i = C_i$
for $1 \leq i \leq n$.  Let $\stackrel{\circ}{\cal G}=
\displaystyle{\oplus_{\alpha \in \stackrel{\circ}{\triangle} \ U \{0\}}}
\stackrel{\circ}{\cal G}_{\alpha}$ be the root space decomposition.
Define
$$
\begin{array}{lll}
\tau_{\alpha+ \delta_{\underline{m}}} &=& \stackrel{\circ}{\cal
G}_{\alpha} \otimes t^{\underline{m}}, \alpha \in
\stackrel{\circ}{\triangle}
\\
 \tau_{\delta_{\underline{m}}} &=&
\stackrel{\circ}{h} \otimes t^{\underline{m}} \oplus \sum_{i=1}^{n}
\C t^{\underline{m}} K_i, \underline{m} \neq 0 \\
\tau_0 &=& \underline{h}
\end{array}
$$
Then clearly $\tau = \displaystyle{\oplus_{\alpha \in \triangle}}
\tau_{\alpha}$ is the root space decomposition with respect to
$\underline{h}$ and consistent with the root system $\triangle$ defined in
section 1.   For $\alpha \in \stackrel{\circ}{\triangle}$ let $X_{\alpha}$ be
the root vector of root $\alpha$.

\paragraph*{Definiton}  A module $V$ of $\tau$ is called integrable
if \\

1) (weight module) $V= \displaystyle{\oplus_{\lambda \in
\underline{h}^*}} V_{\lambda}$
$$ {\rm where} \ V_{\lambda} = \{v \in V \mid h v = \lambda (h) v, \forall
h \in \underline{h} \}$$

2) for all $\alpha \in \stackrel{\circ}{\triangle}, \underline{m} \in \Z^n$
and $v \in V$ there exists an integer $k= k (\alpha, m, v)$ such that
$(X_\alpha (\underline{m}))^k v =0$.

For a weight module $V$ let
$P(V)$ denote the set of all weights.  The following is very standard.

\paragraph*{(2.3) Lemma}  Let $V$ be irreducible integrable module
for $\tau$ with finite dimensional weight spaces.
\begin{enumerate}
\item[(1)] $P(V)$ is $W-$invariant
\item[(2)] dim $V_{\lambda} = dim V_{w \lambda} $ for $w \in W, \
\lambda \in P(V)$
\item[(3)] $\alpha$ real in $\triangle, \lambda \in  P (V)$  then $\lambda
(\alpha^\vee) \in \Z$
\item[(4)] If $\alpha$ is real, $\lambda \in P (V)$ and $\lambda
(\alpha^\vee)>0$ then $\lambda- \alpha \in P (V)$.
\item[(5)]  $\lambda (C_i)$ is a constant (integer) $\forall \lambda
\in P(V)$.
\end{enumerate}
\paragraph*{Proof}  Let $\alpha \in \stackrel{\circ}{\triangle}$ and let
$\alpha= \displaystyle{\sum_{j=1}^{d}} c_j \alpha_j$.  Define
$t_{\alpha} = \sum c_j \frac{\mid \alpha_j \mid^2}{2} \alpha_j^\vee$ so
that $\alpha^\vee = \frac{2 t_{\alpha}}{\mid \alpha \mid^2}$.  One can
easily check that $(t_{\alpha}, t_{\alpha}) = (\alpha, \alpha)$ and
that $\alpha (h) = (t_{\alpha}, h)$ where $t_{\alpha}$ is
the unique with that property because the form (,) on
$\stackrel{\circ}{h}$ is nondegenerate.  Fix $X_{\alpha} \in
\stackrel{\circ}{\cal G}_{\alpha}$ and choose $Y_{\alpha} \in
\stackrel{\circ}{\cal G}_{- \alpha}$ such that $(X_{\alpha},
Y_{\alpha}) = \frac{2}{(t_{\alpha_1} t_{\alpha})}$.

\paragraph*{Claim} $[X_{\alpha}, Y_{\alpha}]= \alpha^\vee$. Consider $(h,
[X_{\alpha}, Y_{\alpha}]) = \alpha (h) (X_{\alpha}, Y_{\alpha})=
\alpha (h) \frac{2}{(t_{\alpha}, t_{\alpha})}$.  By the uniqueness
of $t_\alpha$ the claim follows.

Now $[\alpha^\vee, X_{\alpha} ] = \alpha (\alpha^\vee) X_{\alpha} =
2X_{\alpha}, [\alpha^\vee, Y_{\alpha}]= -2 Y_{\alpha}$.
Thus $X_{\alpha}, Y_{\alpha}, \alpha^\vee$ is an $sl_2$ copy.  Consider
$\gamma= \alpha+\delta_{\underline{m}}$ and recall the definition of
$\gamma^\vee = \alpha^\vee + \frac{2}{\mid \alpha \mid^2} \sum {m_j} C_j$.
Consider $[X_{\alpha} \otimes t^{\underline{m}}, Y_{\alpha} \otimes
t^{-m}]$,
$$
\begin{array}{lll}
&=&\alpha^\vee+ (X_{\alpha}, Y_\alpha) \sum m_i K_i \\
&=& \alpha^\vee+ \frac{2}{\mid \alpha \mid^2} \sum m_i K_i = \gamma^\vee
\end{array}
$$
$$[\gamma^\vee, X_{\alpha} \otimes t^{\underline{m}} ] = 2 X_{\alpha}
\otimes t^{\underline{m}} $$
$$[\gamma^\vee, Y_{\alpha} \otimes t^{- \underline{m}}] = -2 Y_{\alpha}
\otimes t^{- \underline{m}}$$

Thus $X_{\alpha} \otimes t^{\underline{m}}, Y_{\alpha} \otimes t^{-
\underline{m}}, \gamma^\vee$ is an affine $sl_2$.

Now $2,3$ and $4$ follows from standard $sl_2-$theory
representation.  Since $C_i$ is integer linear  combination of
$\alpha^\vee_j$ it follows that $\lambda (C_i)$ is an integer.  The fact
that $\lambda (C_i)$ is constant follows from Lemma 4.3 (2). [QED].

Consider the affine Lie-algebra ${\cal G}_{af}=
\stackrel{\circ}{\cal G} \otimes \C [t_1, t_{1}^{-1}] \oplus \C C_1
\oplus \C d_1$.  Let ${\cal G}_{af}= N^- \oplus h_{1}^{1} \oplus
N^+$ be  the standard decomposition into positive root spaces,
negative root spaces and a Cartan $h^1_1$ spanned by $\alpha_1^\vee,
\cdots, \alpha_d^\vee, \alpha_{d+1}^{\vee}, d_1$.

\paragraph*{(2.4) ~ Proposition}  Let $V$ be irreducible integrable
module for $\tau$.  Assume that $C_1$ acts by non-zero $K_1$ and
$C_i$, $2 \leq i \leq n$ acts trivially.  Then there exists a weight
vector $v$ in $V$ such that
$$N^+ \otimes A_{n-1} v =0 \ {\rm or} \ N^- \otimes A_{n-1} v =0$$
where $A_{n-1} = \C [t_2^{\pm}, \cdots, t_{n}^{\pm}]$ and
$$ N^+ \otimes A_{n-1} = \bigoplus_{\stackrel{m_1 \geq 0 \ {\rm if} \
\alpha \in \stackrel{\circ}{\triangle}^+}{{\rm and} \ m_1>0 \ {\rm
if} \  \alpha \in \stackrel{\circ}{\triangle}^- \cup \{0 \}}}
{\cal G}_\alpha \otimes \ {t_1^{m_1} t^{\underline{m}}}$$

$$\underline{m} = (m_2, \cdots, m_n) \in \Z^{n-1}$$
We need some notation and few lemmas for this.  Let $\lambda \in
P(V)$ then $\lambda$ can be uniquely written as
$$
\begin{array}{lll}
\lambda &=& \overline{\lambda} + \sum  {g}_i \delta_i + \sum s_i
w_i \\
{\rm where} \ \overline{\lambda}& =& \sum_{i=1}^{d} m_i \alpha_i \in
\stackrel{\circ}{h}^*.  \\
{\rm Now \  for} \ 2 \leq i \leq n.&& \\
0 &=& \lambda (C_i) =s_i w_i (C_i) = s_i. \\
K_1 &=& \lambda (C_1) = s_1 w_1 (C_1) = s_1
\end{array}
$$
\paragraph*{(2.5)} Thus $\lambda = \overline{\lambda} + \sum g_i \delta_i
+ K_1 w_1$.

Let $\Gamma_0^+$ be  non-negative integral linear
combination of $\alpha_1, \cdots, \alpha_d$.

\paragraph*{(2.6) ~ Lemma}  Let $V$ be a module satisfying the
conditions of the above proposition.  Then there exists $\lambda \in
P (V)$ such that $\lambda+ \eta \notin P(V), \forall \eta \in
\Gamma_0^+ \backslash \{0\}$ where $(\lambda+ \eta) (d_i) = \lambda
(d_i) = \lambda_1
(d_i)$  for fixed $\lambda_1 \in P(V)$.

\paragraph*{Proof}  Let $\lambda_1 \in P (V)$, let $\lambda_1 (d_i) =
g_i$ let $\underline{g} = (g_1, \cdots g_n)$ and let
$$V_{\underline{g}} = \{v \in V; d_i v = g_i v \}.$$
Let $P_{\underline{g}} (V) = \{ \lambda \in P (V); V_{\lambda} \subseteq
V_{\underline{g}} \}$.  Then from (2.5) for $\lambda \in
P_{\underline{g}} (V)$ we have
$$\lambda = \overline{\lambda} + \sum g_i \delta_i + K_1 w_1.$$
Further $\lambda \mid \stackrel{\circ}{h} = \overline{\lambda}$.  So
if $\lambda, \mu \in P_{\underline{g}} (V)$ such that $\lambda \neq
\mu$ then $\overline{\lambda} \neq \overline{\mu}$.  Now
$V_{\underline{g}}$ is an integrable $\stackrel{\circ}{\cal
G}$-module with finite dimensional weight spaces (with respect to
$\stackrel{\circ}{h})$.  Hence
$$V_{\underline{g}} = \oplus_{\overline{\lambda} \in
\stackrel{\circ}{h}^*} V(\overline{\lambda}) \ {\rm where } \ V
(\overline{\lambda})$$
is an irreducible finite dimensional module for
$\stackrel{\circ}{\cal G}$.  For $\lambda, \mu \in P_{\underline{g}}
(V)$ we have $\lambda - \mu = \displaystyle{\sum_{i=1}^{d}} m_i
\alpha_i$.  Since $V$ is irreducible it follows that $m_i$'s are
integers.  So that $\lambda$'s that occur in $P_{\underline{g}} (V)$
determine a unique coset in the weight lattice of
$\stackrel{\circ}{\cal G}$ modulo its root lattice.  Let $\lambda_0
\in \stackrel{\circ}{h}^*$ be the unique miniscule weight (see Ex 13 of
chapter III of [H]).  Then by Lemma B, 13.4 of [H] it follows then
$\lambda_0$ is a weight of $V (\overline{\lambda})$ that occur in
$V_{\underline{g}}$.  Since each weight space is finite dimensional
the number of $V(\overline{\lambda})$  that occur in
$V_{\underline{g}}$ has to be finite.  Thus $P_{\underline{g}} (V)$
is finite.  Let $\lambda$ be maximal with the ordering
$\frac{<}{\circ}$.
 Then $\lambda + \eta \notin P(V)$ for any $\eta \in \Gamma_+^0
\backslash \{0\}$.  Just note that if $\lambda + \eta \in P (V)$ then
$\lambda+ \eta \in P_{\underline{g}} (V)$.  Further $(\lambda + \eta)
(d_i) = g_i= \lambda (d_i)$. [QED].

\paragraph*{(2.7) } Now we will define a different ordering $\leq $
on $\underline{h}^*$  by  $\lambda \leq \mu$ if $\mu - \lambda =
\displaystyle{\sum_{i=1}^{d+1}} n_i \alpha_i, n_i \in \N$.  Thus if
we say $\lambda \geq 0$ then $\lambda =
\displaystyle{\sum_{i=1}^{d+1}} n_i \alpha_i, n_i \in \N$.

\paragraph*{(2.8) Lemma}  Let $V$ be as in the above proposition with
additional assumptions that $K_1 >0$.  Suppose for all $\lambda \in P
(V)$ there exists $0 \leq \eta \neq 0$ such that $\lambda+ \eta \in
P(V)$.  Then there exists infinitely many $\lambda_i \in P (V), i \in
\Z_+$ such that
\begin{enumerate}
\item[(1)] There exists a wieght vector $u_i$ of weight $\lambda_i$
such that $N^+ u_i=0$.
\item[(2)] $ \lambda_i (d_j) =  \lambda_{k} (d_j) \ \forall i, k \in
\Z_+ \ {\rm and} \ 2 \leq j \leq n $

\item[(3)] $\lambda_i (d_1) -  \lambda_{i-1} (d_1) \in \Z_+$ for $i
\geq 2$.
\item[(4)] There exists a common weight in $V (\lambda_i) \  \forall i$
where $V (\lambda_i)$ is the irreducible highest module generated
by $u_i$.
\end{enumerate}

\paragraph*{Proof}  Let $\lambda$ be an in Lemma 2.6.  Note that
$\lambda (\alpha_i^\vee) \in \N$ for $i=1,2, \cdots, d$.  Follows from
Lemma 2.3(3). Recall that $\stackrel{\circ}{\triangle}$ is finite
root system.  Let $\triangle^{+ a}_{r e} = \{\alpha + m \delta_1,
\alpha \in  \stackrel{\circ}{\triangle}, n \geq 0$ and $ \alpha \in
\stackrel{\circ }{\triangle}^-, n >0 \}$.  Let $\triangle  (\lambda)
= \{ \gamma \in
\triangle^{+ a}_{r e}; \lambda (\gamma^\vee) \leq 0 \}$.  Since $\lambda
(C_1) = K_1 > 0$ it is easy to see that $\triangle (\gamma)$ is
finite.   This situation is very similar to the proof of Theorem 2.4
(i) of [C].  (Here we need our Lemma 2.6 as the arguments of Lemma
2.6 of [C] are not correct.  The 6th line from above on page 322 does
not follows from earlier argument).  Now from the claims 1,2 and 3 in
the proof of 2.4 (i) of [C] we have a vector $v $ in  $V_{\lambda+
p_1 \delta_1},  p_1 \geq 0$ such that
$$\tau_{r \delta_1} v =0, r >0$$
$$\tau_{\alpha+ s \delta_1 } v =0$$

for all but finitely many roots $\alpha+ s \delta_1$ in
$\triangle^{+a}_{r e}$.  Since $V$ is integrable it follows that
$U (N^+) v$ is finite dimensional.  Choose $\lambda_1 $ maximal in
the ordering $\leq$ among the weights in $U (N^+)v$.  Then there
exists a vector $u_1 $ of weight $\lambda_1$ such that
$$N^+ u_1 =0.$$
Further $\lambda_1 = \lambda+ p_1 \delta_1+ \eta_1, \eta_1 \geq 0$
and $\lambda_1 (d_j) = \lambda (d_j)$ for $2 \leq j \leq n$ and
$$\lambda_1 (d_1) - \lambda (d_1) = p_1+\eta_1 (d_1) \geq 0.$$
Now by Lemma 2.6 there exists $\tilde{\lambda}_1 \in P (V)$ such that

\paragraph*{(2.9)} $ V_{\tilde{\lambda}_1+ \eta} =0 ~~ \forall \eta
\in \Gamma_0^+ \backslash \{0\}$ and $ \tilde{\lambda}_1 (d_i) =
\lambda_1 (d_i) \ 1 \leq i \leq n$.  By the assumption in the Lemma
there exists $\tilde{\eta}_1 >0$ such that
$$\tilde{\lambda}_1 + \tilde{\eta}_1 \in P (V).$$
By (2.9) it follows that $\tilde{\eta}_1 (d_1) > 0$. As $
\tilde{\eta}_1 \geq 0$ and we also have $\tilde{\eta}_1 (d_i) = 0$
for $2 \leq i \leq n$.  Repeating the above argument for
$\tilde{\lambda}_1+ \tilde{\eta}_1$ in place of $\lambda$ to get a
weight vector $u_2$ of weight $\lambda_2 =
\tilde{\lambda}_1+\tilde{\eta}_1 +p_2 \delta_1+ \eta_2, \eta_2
\geq 0$ such that $N^+ u_2 =0$.  Now $\lambda_2 (d_j) =
\tilde{\lambda}_1 (d_j) = \lambda (d_j)$ for $ 2 \leq j \leq n$.
$$\lambda_2 (d_1) = \lambda_1 (d_1) + p_2 + \eta_2 (d_1) +
\tilde{\eta}_1 (d_1)$$ so that $\lambda_2 (d_1) - \lambda_1 (d_1)
= p_2 +\eta_2 (d_1)+ \tilde{\eta}_1 (d_1) >0$. (Note the strict
inequality).  Repeating this process we have (1), (2) and (3).
Clearly $\lambda_i ({\alpha}_j^\vee) \in \N \ \ 1 \leq j \leq
1+d$.  By irreducibility of $V$ we have
$$\overline{\lambda}_i - \overline{\lambda}_j = \sum_{i=1}^{d} m_i
\alpha_i \ {\rm for} \ m_i \in \Z.$$
Then $\overline{\lambda}_i$s determine unique coset of weight
lattice.  As earlier let $\overline{\lambda}_0 \in
\stackrel{\circ}{h}^*$ be the unique miniscule weight.  Thus
$\overline{\lambda}_0 \frac{<}{0} \overline{\lambda}_i \forall i \in
\Z$.  (Note that each $\overline{\lambda}_i$ is dominant integral
weight).

Let $\lambda_0 = \overline{\lambda}_0 + \displaystyle{\sum_{i=1}^{d}}
\lambda (d_i) \delta_i +K_1 w_1 \leq \lambda_i$.  This
inequality is by construction of $\lambda_i$.  Recall $\lambda_i=
\overline{\lambda}_i +\lambda_i (d_1) \delta_1+
\displaystyle{\sum_{i=2}^{d}} \lambda (d_i) \delta_i + K_1 w_1$
and we have $\lambda_i (d_i) - \lambda (d_1)>0$.  Now $\lambda_0
(\alpha_i^\vee) = \overline{\lambda}_0 (\alpha_i^\vee) \in \N$ for $ i=1,2,
\cdots,d$.
$$
\begin{array}{lll}
\lambda_0 (\alpha_{d+1}^{\vee}) &=& \lambda_0 (C_1 - \beta^\vee) \\
&=& K_1 - \overline{\lambda}_0 (\beta^\vee)
\end{array}
$$
Since $K_1$ is integer and positive we have $K_1 \geq 1$.  By Lemma
(1.2) we have $\overline{\lambda}_0 (\beta^\vee)=0$ on 1.  Thus
$\lambda_0 (\alpha_{d+1}^{\vee}) \in \N$.  Hence $\lambda_0$ is dominant
integral and $\lambda_0 \leq  \lambda_i$.  By Proposition 12.5 (a) of
[K] it  follows that $\lambda_0$ is a weight of $V(\lambda_i)$ for all
$i$.  This proves (4). [QED].

\paragraph*{Proof of the Proposition (2.4)}  Assume $K_1 >0$.
Suppose the conclusions of Lemma (2.8) are true.  Since $\lambda_i$'s
are distinct, $V(\lambda_i)$ are all non-isomorphic irreducible
heighest weight modules.  Hence their sum has to be direct.  But
$V_{\lambda_0} \cap V (\lambda_i) \neq 0$.  Hence dim $V_{\lambda_0}$
is infinite.  A contradiction.

So we conclude that there exists a $\lambda \in P (V)$ such that
\paragraph{(2.10)} $ \lambda+ \eta \notin P(V), \forall {0 \neq}  \eta
\geq 0$.  In particular $\lambda+ \alpha \notin P(V)$ for $\alpha \in
\triangle^{+ a}_{re}$.  Thus from Lemma (2.3) (4) it follows that
$\lambda (\alpha^\vee) \in \N$:

\paragraph*{(2.11)} Suppose $V_{\lambda + \alpha + \delta'}=0$ for
all $\alpha \in \triangle^{+ a}_{re}$ and for all $\delta^1 =
\displaystyle{\sum_{i=2}^{n}}  m_i \alpha_i$.  This means $X_{\alpha}
\otimes  t^{\underline{m}} V_{\lambda} =0$.  But the positive real
roots generates the positive null roots and hence $N^+ \otimes A_{n-1}
V_{\lambda} =0$.  So we are done.  Now assume that $V_{\lambda+
\alpha + \delta} \neq 0$ for some $\alpha \in \triangle^{+ a}_{re}$
and for some $\delta= \displaystyle{\sum^{n}_{i=2}} n_i \delta_i$.
Let $\mu= \lambda+ \alpha+ \delta$.

\paragraph*{Claim} $ V_{\mu+ \gamma+ \delta^1}= 0$ for $ \gamma \in
\triangle^{+a}_{re} $ and for all $\delta^1=
\displaystyle{\sum_{i=2}^{n}} m_i \delta_i$.  Suppose it is false.
Then $V_{\mu+ \gamma+ \delta^1}  \neq0$ for some $\gamma$ and
$\delta^1$.

\paragraph*{Case 1}  $(\alpha+ \gamma, \alpha) >0$.  We know by Lemma
(1.4) (2) that $(\alpha+ \gamma)(\alpha^\vee) >0$. Put $\gamma_1 =
\alpha+ \delta+ \delta^1$ and consider
$$(\mu+ \gamma+ \delta^1) (\gamma_1^\vee) = (\lambda+ \alpha+ \gamma)
(\alpha^\vee)>0.$$
We are using the fact that $\lambda (C_i) =0$ for $2 \leq i \leq n$.
By Lemma (2.3) (4), $\lambda+\alpha+ \gamma +\delta+\delta^1 - (\alpha+
\delta+\delta^1) = \lambda+\gamma \in P (V)$ a contradiction to
(2.10).

\paragraph*{Case 2}  $(\alpha+\gamma, \gamma)>0$ which can be done in
the same way.

Since the symmetric billinear form given by an affine
matrix is semi positive definite on the root lattice, we are left
with the

\paragraph*{Case 3} $(\alpha+\gamma, \alpha+\gamma)=0$.  This implies
$\alpha+\gamma= \ell \delta_1$ for $\ell >0$ and further
$(\alpha+\gamma, \alpha) =0$ and $(\alpha+\gamma, \gamma)=0$.  Then
$(\lambda+ \alpha+ \gamma +\delta+\delta^1) (\gamma_1^\vee) = \lambda
(\alpha^\vee)$ where $\gamma_1= \alpha+\delta+ \delta^1$.

\paragraph*{Sub case 1} $\lambda (\alpha^\vee)>0$ (in any case $\lambda
(\alpha^\vee) \geq 0)$.  By lemma 2.3 (4) we have
$$\lambda+ \alpha+ \gamma+ \delta+\delta^1 - (
\alpha+\delta+\delta^1) = \lambda+ \gamma \in  P (V)$$
a contradiction to (2.10).

\paragraph*{Sub case 2}  $\lambda (\alpha^\vee)= 0$.
Note that $\lambda((- \alpha+ \ell \delta_1 + \delta+ \delta^1)^\vee )$
$$= \ell \frac{2}{\mid \alpha \mid^2} \lambda (C_1) >0$$
Hence by Lemma 2.3 (4) we have
$$\lambda+ \alpha+\gamma + \delta+\delta^1 -
(\gamma+\delta+\delta^1)=\lambda+ \alpha \in P (V)$$
a contradiction to (2.10).  Thus our claim follows.  As earlier it
follows that $N^+ \otimes A_{n-1} V_{\mu}=0$.  Hence we are done.
The case $K_1<0$ can be done similarly.[QED].

\paragraph*{(2.12) Proposition} Let $V$ be integrable irreducible
module for $\tau$ with finite dimensional weight spaces.  Suppose
$K_i=0$ for $ 1 \leq i \leq n$.  Then there exists a weight vector
$v$ in $V$ such that
$$n^+\otimes A_n v=0 \ {\rm or} $$
$$n^- \otimes A_n v=0$$
where $\stackrel{\circ}{\cal G}=n^+ \oplus \stackrel{\circ}{h} \oplus
n^-$ is the standard decomposition.

\paragraph*{Proof}  This follows from the proof of theorem 2.4 (ii)
of [C].  Use Lemma (2.6) of our paper instead of Lemma 2.6 of [C].
The fact that the dimension of null roots could be greater than one
does not matter.

\section*{3 Section}  Let $ {\cal G}$ be a Kac-Moody Lie-algebra.
Let $h$ be a Cartan subalgebra of ${\cal G}$.  Let ${\cal G}^1=
[{\cal G}, {\cal G}]$ and let $h^1= h \cap {\cal G}^1$.  Let ${\cal
G}^1 = N^+ \oplus h^1 \oplus N^-$ be the standard decomposition into
positive roots spaces, negative roots spaces and a $h^1$.  Fix a
positive integer $n$ and let $A=A_n = \C [t_1^{\pm}, \cdots,
t^{\pm}_{n}]$.  Let $D$ be a linear span of derivations $d_1
\cdots, d_n$.  Then let $\tilde{{\cal G}}_A= {\cal G}^1_A \oplus h^{''}
\oplus D$ where $h^{"}$ is define as $h=h' \oplus h^{"} $ (see $\S$
1.3 of [K]).  Define Lie-algebra structure on $\tilde{\cal G}_A$ by

\begin{eqnarray*}
[X (\underline{r}), Y (\underline{s})] &=& [X,Y]
(\underline{r}+\underline{s}). \\ %
\left [ d_i, X(\underline{r}) \right ] &=& r_i X (\underline{r}), X,Y
\in {\cal G} \
\underline{r}, \underline{s} \in \Z^n \\ %
\left [h, X (\underline{r})\right ]&=& [h, X] (\underline{r}) \\ %
\left [h, d_i \right ] &=& [d_i,d_j]=0, h \in h
\end{eqnarray*}
Let $\tilde{h}_A=h^1 \otimes A \oplus D \oplus h^{''}$ be an abelian
Lie-algebra.  For any Lie-algebra $G$, let $U(G)$ be the universal
enveloping algebra.  Let $H=h^1 \otimes A \oplus h^{''}$.  Then $U(H)$ is
clearly $\Z^n-$graded abelian-Lie-algebra.  Let $\overline{\psi}:U
(H) \to A $ a $\Z^n-$graded homomorphism.  Then $A$ is a module for
$H$ via $\overline{\psi}$ defined as
$$
\begin{array}{lll}
h (\underline{m}) t^{\underline{s}} &=& \overline{\psi} (h
(\underline{m})) t^{\underline{s}} \  \ h \in H \ \underline{m} \in \Z^n \\
h t^{\underline{s}} &=& \overline{\psi} (h) t^s \ {\rm for} \ h \in h^{''}
\end{array}
$$
Let $A_{\overline{\psi}} =$ Image of $\overline{\psi}$.

\paragraph*{(3.1) Lemma} (Lemma (1.2) , [E1]).$ A_{\overline{\psi}}$ is
an irreducible $\tilde{h}_A-$module if and only if each homogeneous
element of $A_{\overline{\psi}}$ is invertible.

Just note that $h^{''}$ does not play any role.

We need the notion of highest weight module for $\tilde{\cal G}_A$.  Let
$\overline{\psi}$ be as above.

\paragraph*{(3.2) Definition}  A module $V$ of $\tilde{\cal G}_A$ is
called highest weight module for $\tilde{\cal G}_A$ if there exist a
weight vector (with respect to $h \oplus  D) \  v$ in $V$ such that

(1) $V = U (\tilde{\cal G}_A)v$ \ (2) $N_A^+ v=0$ \ (3) $U
(\tilde{h}_A)v$ is an irreducible module for $\tilde{h}_A$
given by $\overline{\psi}$.

We will now define universal highest weight module for $\tilde{\cal
G}_A$.  Let $\overline{\psi}$ be as above such that
$A_{\overline{\psi}}$ is irreducible.  Let $N_A^+$ act trivially on
$A_{\overline{\psi}} $.  Now consider the following induced
$\tilde{{\cal G}}_A-$ module.
$M (\overline{\psi}) = U (\tilde{{\cal G}}_A) \otimes_B
A_{\overline{\psi}}$ where  $B= N^+_A \oplus \tilde{h}_A$.

We have the following standard.

\paragraph*{(3.3) Proposition}
\begin{enumerate}
\item[(1)] As $h \oplus D$ module, $M (\overline{\psi})$ is a weight
module.
\item[(2)] $M (\overline{\psi})$ is a free $N^-_A$ module and as a
vector space
$$M (\overline{\psi}) \cong U (N^-_A) \otimes A_{\overline{\psi}}.$$
\item[(3)]  $M (\overline{\psi})$ has a unique irreducible quotient
called $V(\overline{\psi})$.
\end{enumerate}
\paragraph*{Proof}  See [E3].

We need the following non-graded highest weight module for ${\cal
G}_A^1 \oplus h^{''}$.

\paragraph*{(3.4) Definition}  A module $W$ of ${\cal G}^1_A \oplus
h^{''}$ is said to be (non-graded) highest weight module if there
exists a weight vector (with respect to $h) \ v$ in $W$ such that
(1) $U ({\cal G}^1_A \oplus h^{''}) \ v=W$ \\
(2) $N_A^+ v =0$ \\
(3) There exists a $\psi$ in $H^*$ such that $h v= \psi (h) v$ for
all $h \in H$.

Let $\psi \in H^* $ and let $H$ act as one dimensional vector space
$ \C (\psi)$ by $\psi$.  Let $N^+_A$ acts trivially on $\C (\psi)$.
Consider the induced module $M (\psi) = U ({\cal G}^1_A \oplus h^{''})
\displaystyle{\otimes_{B'}} \C (\psi)$ where $B'= N^+_A \oplus H$.

Now by standard arguments $M (\psi)$ has a unique irreducible
quotient $V(\psi)$.

Let $\overline{\psi}$ be as above and let $\psi= E(1) \circ
\overline{\psi}$ where $E(1): A_{\overline{\psi}} \to \C$ defined by
$E (1) t^{\underline{m}} = 1$.  We will make $V(\psi)_A$ a (graded) $
\tilde{{\cal G}}_A-$module by \\
$g (\underline{m}) v (\underline{r}) = (g (\underline{m}) v)
(\underline{m}+\underline{r})$ for $g \in {\cal G}^1_A,
\underline{r}, \underline{m} \in \Z^n, v \in V (\psi)$. \\
$$d_i v (\underline{r}) = r_i v (\underline{r})$$
$$h^{''} v (\underline{r})= (h^{''} v) (\underline{r}), h^{''} \in
h^{''}.$$

\paragraph*{(3.5) Proposition}  Let $\overline{\psi}$ and $\psi$ as
above.  Assume that $A_{\overline{\psi}}$ is irreducible
$\tilde{h}_{A^-}$ module.  Let $G \subset \Z^n$ be such that
$\{t^{\underline{m}}, \underline{m} \in G\}$ is a coset
representatives for $A/A_{\overline{\psi}}$.  Let $v$ be a highest
weight vector of $V(\psi)$.  Then \\
$ (1) ~~~~ V (\psi)_A = \oplus_{\underline{m} \in G} U v
(\underline{m})$ as  $\tilde{{\cal G}}_A-$module.  $U v
(\underline{m})$ is a  $\tilde{{\cal G}}_A-$module generated by $v
(\underline{m})$.

(2) Each $Uv (\underline{m})$ is an irreducible $\tilde{{\cal
G}}_A-$module.

(3) $Uv (0) \cong V (\overline{\psi})$ as $\tilde{{\cal
G}}_A-$module.

\paragraph*{Proof}  Follows from Theorem (1.8) of [E1].  It is stated
for special $\overline{\psi}$.  But we have only used the fact
$A_{\overline{\psi}}$ is irreducible  $\tilde{h}_{A^-}$module.  More
over our $\tilde{{\cal G}}_A$ is smaller than one in [E] but all the
arguments take place inside $\tilde{\cal G}_A$. [QED].

\paragraph*{(3.6) Lemma}  (Lemma (1.10) , [E3])  Let $\psi$ and
$\overline{\psi}$ as above.  Then $V(\overline{\psi})$ has finite
dimensional weight spaces (with respect to $h \oplus D)$ if and only
if $V(\psi)$ has finite dimensional weight spaces (with respect to
$h$).

The following Lemma tells for which $\psi, V(\psi)$ has finite
dimensional weight space with respect to $h$.  There by giving
conditions for $V (\overline{\psi})$ to have  finite dimensional
weight spaces.

\paragraph*{(3.7) Lemma}   $V (\psi)$ has finite dimensional weight
spaces if and only if $\psi$ factors through $h' \otimes A/I$ for some
co-finite ideal $I$ of $A$.

\paragraph*{Proof}  Assume that $V(\psi)$ has finite dimensional
weight spaces.  Let $\triangle^+$ be a positive root system for
${\cal G}$.  Let $\alpha$ be a  simple root in $\triangle^+$.  Let
$Y_{\alpha}$ be a root vector for the root $- \alpha$.  For fixed $j,
$ consider, for a highest weight vector $v$,
$$\{Y_{\alpha} \otimes t_{j}^{n} v, n \in \Z \}$$
which is contained in the same weight space $V_{\lambda - \alpha}
(\psi), \lambda = \psi \mid h$.  Thus there exists a polynomial $P_j
(\alpha) = \sum a_i t^i_j$ such that
$$Y_{\alpha} \otimes P_j (\alpha) v =0.$$
Where $Y_{\alpha} \otimes P_j (\alpha) = \sum a_i Y_{\alpha} \otimes
t^i_j$.  Let $(P)$ be the  ideal generated by polynomial $P$ inside $A$.

\paragraph*{Claim 1}  $Y_{\alpha} \otimes (P_j (\alpha)) v=0$. \\
Consider $0 = h (\underline{m}) Y_{\alpha}  \otimes P_j (\alpha) v =
Y_{\alpha} \otimes P_j (\alpha) h (\underline{m}) v
- \alpha (h) Y_{\alpha} \otimes  t^{\underline{m}} P_j (\alpha) v$.
Since $h (\underline{m})$ acts by scalar on $v$, the  first term is
zero.  Hence the claim.

Put $P_j =  \prod P_j (\alpha)$ where $\alpha$ runs through all simple
positive roots.

\paragraph*{Claim 2}  $Y_{\beta} \otimes (P_j) v =0$ for $\beta \in
\triangle^+$.  First note that  $(P_j) \subset P_j (\alpha)$ and
hence the claim is true for any simple root. Claim by induction on the
height $\beta$.  Caim is true for $\beta$ such that  height $\beta=1$.
Let $\alpha$ be a simple root in $\triangle^+$.  Let $X_{\alpha}$ be
the root vector of root $\alpha$.  Consider
$$X_{\alpha} (\underline{m}) Y_{\beta} \otimes (P_j) v = Y_{\beta}
\otimes (P_j) X_{\alpha} (\underline{m}) v+ [X_{\alpha}, Y_{\beta}]
\otimes t^{\underline{m}} (P_j) v.$$
The first term is zero as $v$ is highest weight vector.  The second
term is zero by induction.  Since $Y_{\beta} \otimes (P_j) v$ is
killed by $X_{\alpha} (\underline{m})$ for $\alpha$ simple and for
any $\underline{m} \in \Z^n$ it is easy to see that $X_{\alpha}
(\underline{m}). Y_{\beta} \otimes (P_j) v=0$ for all $\alpha \in
\triangle^+$ and $\underline{m} \in \Z^n$.  Hence $Y_{\beta} \otimes
(P_j) v$ is a highest weight of weight $\lambda -\beta$.  But
$V(\psi)$ is an irreducible highest weight module and hence
$Y_{\beta} \otimes (P_j) v=0$.  Hence claim 2.

\paragraph*{Claim 3}  $h \otimes (P_j) v =0 \ \ \forall h \in h'$. \\
Consider $h_{\alpha} \otimes (P_j) v = X_{\alpha} Y_{\alpha} \otimes
(P_j) v - Y_{\alpha} \otimes (P_j) X_{\alpha} v =0$ \  since
$h_{\alpha}$ covers all $h'$ for $\alpha$ simple  we have claim 3.

Let $I$ be an ideal generated by $P_1, P_2, \cdots, P_n$ inside $A$.
It is elementary to see that $A/I$ is finite dimensional.  Thus
$\psi$ factors through $h' \otimes A / I$.  In fact we can prove that
${\cal G}' \otimes I. V (\psi) =0$ by considering $W= \{w \in V
(\psi), {\cal G}' \otimes  I w=0 \}$.  We have just seen that $W$ is
non empty.  It is easy to see that $W$ is a sub module.  Hence $W= V
(\psi)$.

Conversely suppose that $\psi$ factors through $h' \otimes A/I$ for a
co-finite ideal  $I$ of $A$.

\paragraph*{Claim 4}  For $\beta > 0, Y_{\beta} \otimes I v =0$.  Let
$\alpha$ be simple positive.  Consider
$$
\begin{array}{lll}
X_{\alpha}(\underline{m}) Y_{\alpha} \otimes I v &=& Y_{\alpha}
\otimes I X_{\alpha} (\underline{m})v
+ h_{\alpha} \otimes I v \\
&=& 0
\end{array}
$$
Let $\alpha_1$ be simple positive root different from $\alpha$.
Then clearly $X_{\alpha}(\underline{m})Y_{\alpha} \otimes I v=0$.
Hence $Y_{\alpha} \otimes I v$ is a highest weight vector.  Since
$V(\psi)$ is irreducible highest weight module we conclude that
$Y_{\alpha} \otimes I v =0$.

Now arguing as earlier on the induction of ht $\beta$ we have
$Y_{\beta} \otimes Iv=0$.  Hence the claim.  Thus we have ${\cal G}^1
\otimes I v=0$.  Now consider the non-zero submodule
$$W = \{w \in V(\psi), {\cal G}^1 \otimes I w =0 \}$$
of $V (\psi)$.  Since $V(\psi)$ is irreducible we have $W=V(\psi)$.

Let $V_{\mu} (\psi)$ be a wiehgt space of $V(\psi)$.  Then by $PBW$
theorem any vector  of $V_{\mu} (\psi)$ is linear combination of the
vector of the form.

\paragraph*{(3.8)}  $Y_{\alpha_1} t^{\underline{m}_1} Y_{\alpha_2}
t^{\underline{m}_2} \cdots Y_{\alpha_k} t^{\underline{m}_k}$ such that $- \sum \alpha_i+ \lambda = \mu$
where $\lambda = \psi \mid h$ and $\alpha_i > 0$.  Thus the number of
$\alpha_i$ that can occur are finite.  Now for any finite dimensional
space $W$ of $V(\psi)$ and for any fixed $\alpha_i$ the space
$Y_{\alpha_i} \otimes t^{\underline{m}} W$, for $\underline{m} \in
\Z^n $ is finite dimensional.  Thus at every stage in the equation
(3.8) we get a finite dimensional  space.  So we conclude that $V_{\mu}
(\psi)$ is finite dimensional. [QED].

\paragraph *{(3.9) Remark  (1)} $V (\psi)$ has finite dimensional
weight space if and only if ${\cal G}' \otimes I. V (\psi) =0$ for a
co-finite ideal $I$ of $A$.

(2)  This lemma also gives new modules with finite dimensional weight
spaces for affine Lie-algebras by taking ${\cal G}$ to be finite
dimensional simple Lie-algebra and $n=1$.

We will now give a continious family of irreducible highest weight
modules for ${\tilde{\cal G}}_A$ with finite dimensional weight
spaces.

As earlier let $n$ be a positive integer.  For each $i,  1 \leq i \leq
n$, let $N_i$ be a positive integer.  Let $\underline{a}_i = (a_{i
1}, \cdots, a_{i N_i})$ be non-zero distinct complex numbers.  Let $N=
N_1 \cdots N_n$.  Let $I = (i_1, \cdots, i_n)$ where $1 \leq i_j
\leq N_j$.  Let $\underline{m} =  (m_1, \cdots, m_n) \in \Z^n$.
Define $a_I^{\underline{m}} = a_{1i_1}^{m_1} \cdots a_{ni_n}^{m_n}$.
Let $\phi$ be a Lie-algebra homomorphism defined by

\paragraph*{(3.10) } $ \Phi: {\cal G}_A \to {\displaystyle{\oplus_{N-
{\rm copies} \ }}} {\cal G}= {\cal G}_N$ \\
$$X \otimes t^{\underline{m}} \mapsto (a_I^{\underline{m}} X).$$
${\cal G}$ could be any Lie-algebra. This map was first defined by
Kac-Jacobson for $n=1$.

\paragraph*{Lemma (3.11)}  (a) $\Phi$ is surjective.

(b)  Let $P_j (t_j) = \displaystyle{\prod_{k=1}^{N_j}} (t_j - a_{jk})
\ {\rm and} \ I $ be the ideal generated by $P_1 (t_1), \cdots, P_n
(t_n)$ inside $A$.  Then ${\cal G}\otimes A/I \cong {\cal G}_N$.
\paragraph*{Proof}

(a) We will first prove that the following $N \times N$ matrix is
invertible.
$$
\begin{array}{lll}
X&=& (a_{1i}^{m_1} \cdots a_{ni_n}^{m_n})_{\stackrel{0 \leq m_i \leq
N_i-1}{ 1 \leq i_j \leq N_j.}}
\end{array}
$$

The index $\underline{m}$ determines rows and the index $(i_1, \cdots,
i_n)$ determines columns. For $n=1, X$ becomes Vandermonde matrix as
$a_{1i} \neq a_{1j}$.  Hence $X$ is invertible.  We will now prove
this for $n=2$ and then extend it for all $n$.

Given a square matrix $A$, we call a square matrix of the form.
$$
\begin{pmatrix}
A \ 0 \ \cdots \ 0 \cr 0 \ A \ \cdots \ 0 \cr \cdots \ \cdots \
\cdots  \cdots \cr 0 \ 0 \ \cdots \ 0 A
\end{pmatrix}
$$
a block diagonal matrix of $A$.  Let $B=(\stackrel{j}{a_{1 i}})_{1
\leq i, j \leq N_1} \ {\rm and} $ \\
$$C= (\stackrel{j}{a_{2i}})_{1 \leq i, j \leq N_2}.$$
For $K=1,2$, since $a_{Ki}$ are distinct non-zero complex numbers for
distinct $i, B$ and $C$ are Vandermonde matrices and hence they are
invertible.  Let $\tilde{B}$ be $N \times N$ block diagonal matrix of
$B$.  Similarly $\tilde{C}$ for $C$.  Clearly both $\tilde{B} $ and
$\tilde{C}$ are invertible.

Let $\sigma:<1 ,2, \cdots, N > \to < 1, \cdots, N >$ be such that
$$\sigma (p) = (m-1) N_2 + (\ell+1)$$
where $p =N_1 \ell +m, 0 \leq \ell \leq N_2 -1, 1 \leq m \leq N_1$
(recall that $N_1 N_2 =N)$. Clearly $\sigma$ is injective.
Let $\tilde{D} = ( \stackrel{\approx}{c}_{k, \ell})$ where $
\stackrel{\approx}{c}_{k, \ell} = \tilde{c}_{\sigma (k), \ell},
\tilde{C} = (\tilde{c}_{k, \ell})$ and
$\tilde{B} = (\tilde{b}_{k, \ell})$.  Then clearly $\tilde{D}$ is a
product of permutation matrix and $\tilde{C}$ and hence invertible.
We now claim that upto permutation matrix $X=\tilde{B} \tilde{D}$ which
is invertible.

Consider $(k, \ell)$ entry of $\tilde{B} \tilde{D}$.

\paragraph*{(3.12) } $\displaystyle{\sum_{j=1}^{N}}
\tilde{b}_{kj} \tilde{c}_{\sigma (j) \ell}$

\paragraph*{Claim}  Exactly one and only one term in (3.12) is non zero.
\begin{eqnarray*}
{\rm let} \ k= p_1 N_1+ q_1, && 0 \leq p_1 \leq N_2 -1 \\
&&1 \leq q_1 \leq N_1 \\
\ell = sN_2+ t, && 0 \leq s \leq N_1 -1 \\
&&1 \leq t \leq N_2
\end{eqnarray*}

Since $\tilde{B}$ is a block diagonal matrix of $B, \tilde{b}_{kj}$
to be non zero in (3.11), $j$ should be equal to $p_1 N_1+q_2$ for
some $1 \leq q_2 \leq N_1$.

Now $\sigma (j) = (q_2-1)N_2+(p_1+1)$.  Now since $\tilde{C}$ is a
block diagonal matrix, $\tilde{C}_{\sigma (j), \ell}$ to be non-zero
$\ell$ should be equal to

$(q_2 -1) N_2 +t$ which forces $q_2 -1 =s$.  Thus (3.12) equal to
$$a_{1 q_1}^s a_{2 (p_1+1)}^{t}.$$
This proves the claim.  Notice that the exact placement of the
entries are given by the permutation matrix.  But the entries of the
rows and  columns do not get mixed up.  Hence $X= \tilde{B} \tilde{D}$
up to permutation.

Now we will prove the result for any $n$. Let $A_k =
(a^{j-1}_{ki})_{1 \leq i, j \leq N_k}$ and let $\tilde{A}_k$ be a
block diagonal matrix of $A_k$.  Consider $\tilde{A}_1$ and
$\tilde{A}_2$.  Let $E$ be a permutation matrix given as in the case
$n=2$.  The only difference is that $E$ is actually a block diagonal
matrix of a permutation matrix of order $N_1 N_2$.  Now $\tilde{A}_1
E \tilde{B}_2$ is a block diagonal matrix of a matrix of order $N_1
N_2$.  Now take $\tilde{A}_1 E \tilde{A}_2$ in place of $\tilde{A}_1$
and $\tilde{A}_3$ in place of $\tilde{A}_2$ in the case $n=2$.
Repeat the process to get a  matrix of the form $\tilde{A}_1 E
\tilde{A}_2 F \tilde{A}_3 (F$ is a permutation matrix) which is a
block diagonal matrix of a matrix of order $N_1 N_2 N_3$.  Each entry
of this matrix is of the form
$$ a_{1i_1}^{j_1-1} a_{2 i_2}^{j_2-1} a_{3i_3}^{j_3-1}.$$

Repeating this process we get the desired matrix $X$ which is
invertible.
\vskip 4mm

Let $(X_1, \cdots, X_N) \in {\cal G}_N$. \\
Now $ \Phi (X^{-1} (X_1, \cdots, X_N)^T = (X_1, \cdots, X_N)$.  This
completes the proof of Lemma 3.11 (a).

(b) ${\cal G}\otimes P_j (t_j) \subset ker \Phi$, since the
restriction of $\Phi$ is nothing but evaluation map at the roots of
$P_j (t_j)$.  Thus we have ${\cal G}\otimes I \subseteq ker \Phi$.
Consider the space $T= \{ {\cal G}\otimes t^{\underline{m}}, 0 \leq
m_i < N_i \}$.  Since any ${\cal G}\otimes P$ can be reduced to the
linear combinations of elements of $T$ modulo ${\cal G}\otimes I$,
it is a spanning set.  In (a) we actually proved $\Phi$ is surjective
on $T$.  Further $\Phi$  is injective on $T$ as the corresponding
matrix is invertible.  Thus it follows that ${\cal G}\otimes A/I
\cong {\cal G}_N$. [QED].

\paragraph* {(3.13) Remark}  The condition $a_{ij} \neq a_{ik}$ is
necessary.  Otherwise $\Phi$ is not surjective.

Let $V (\lambda_i)$
be irreducible heighest weight module for ${\cal G}$.  Then it is
known by Lemma (9.10) of [K] that $V (\lambda_i)$ is irreducible
${\cal G}'-$module.  Thus $V= \displaystyle{\otimes_{i=1}^{N}}V
(\lambda_i)$ is an irreducible ${\cal G}_N^{'}-$module.  Restrict the
map $\Phi$ in (3.10) to ${\cal G}'_A \oplus h^{''}$ whose image
containing ${\cal G}'_N$.  Thus $V$ is irreducible ${\cal G}'_A
\oplus h^{''}$ module via $\Phi$.  Consider $I= (i_1, \cdots, i_n)$
for $1 \leq i_j \leq N_j$.  They are $N$ of them.  Give them an order
say $I_1, \cdots I_N$.  The map $\Phi$ is defined in this order.  Now
define \\

\paragraph*{(3.14)} $\psi:U (H) \to \C $ \\
by $ \psi (h \otimes t^{\underline{m}})
= \displaystyle{\sum_{i=1}^{N}} a_{I}^{\underline{m}} \lambda_i (h),
h \in h'$ and $\psi (h) = \sum \lambda_i (h), h \in h$.  Then it is
easy to see that $V$ is a  highest weight
module with highest weight $\psi$.  Thus we have $V (\psi) \cong V$.
Now define

\paragraph*{(3.15)} $\overline{\psi}:U (H) \to A$ as $\Z^n-$graded
homomorphism by \\
$\overline{\psi} (h \otimes t^{\underline{m}}) = \psi (h \otimes
t^{\underline{m}} ) t^{\underline{m}}$.  So in order to get modules
for $\tilde{\cal G}_A$ we need to know that $A_{\overline{\psi}}$  is
an irreducible $\tilde{h}_A-$module.  (See proposition (3.5) and
Lemma (3.1)). It is well known that $V(\lambda_i)$ is  integrable if
and only if $\lambda_i$ is dominant integral.  Since our interest is in
constructing integrable modules we will assume that each $\lambda_i$
is dominant.

\paragraph*{(3.16) Lemma}  Let $\overline{\psi}$ be as above.  Assume
each $\lambda_i$ is dominant and not all of them are zero.  Then
$A_{\overline{\psi}}$ is irreducible $\tilde{h}_A-$module.

\paragraph*{Proof}  Let $I_j= (i_1, \cdots, i_n)$ and let $\lambda_j
= \lambda_{i_1 i_2 \cdots i_n}$.  Consider $\Gamma = \{\underline{m}
\in \Z^n, t^{\underline{m}} \in A_{\overline{\psi}} \} \  \psi$ being
an algebra homomorphism, $\Gamma$ is closed under addition.
$\overline{\psi} (h \otimes 1) = \sum \lambda_{i_1 \cdots i_n}$ and
each weight is dominant and not all zero.  So $\overline{\psi} (h
\otimes 1) \neq 0$ for some $h$.  Hence $0 \in \Gamma$.  To prove the
Lemma it is sufficient to prove that $\Gamma$ is a  sub group (see
Lemma (3.1)).  For $ 1 \leq j \leq N$ consider $\overline{\psi}_j =
\overline{\psi} \mid h \otimes \C [t_j, t_j^{-1}]$.

\paragraph*{Claim}  Image $\overline{\psi}_j = \C [t_j^{k_j},
t_j^{-k_j}]$ for some $k_j \leq N_j$. Consider $\overline{\psi}_j (h
\otimes t_j^{m_j}) = \displaystyle{\sum_{i_j=1}^{N_j}}
\displaystyle{\sum_{(i_1, \cdots, i_n)}} \
\lambda_{i_1, \cdots i_n} (h) a_{ji_j}^{m_j} t_j^{m_j}$.  Clearly for
a fixed $i_j$ one of ${\displaystyle{\sum_{(i_1, \cdots, i_n)}}}
\lambda_{i_1, \cdots, i_n} (h)$ is non zero for some $h$.  Now
$\overline{\psi}_j (h \otimes t_j^{m_j})=0$ for $0 < m_j \leq N_j$
cannot hold as $(a_{ji_j}^{m_j})$ is a Vandermonde matrix.  The same
holds for $-N_j \leq m_j < 0$.  Thus the image has to be
$\C [t_j^{k_j}, t_j^{-k_j}]$ for some  $0 < k_j \leq N_j$.  In fact
$k_j \mid N_j$ see Lemma (5.3) of $[E2]$.  Thus we have
\paragraph*{(3.17)}
$$\ell_1 k_1 \Z \oplus \cdots \oplus \ell_n k_n \Z \subseteq \Gamma$$
for any $\ell_1, \cdots, \ell_n \in \Z$.  $0 < k_i \leq N_i$.  Let
$\Gamma_0 = \{(m_1, \cdots, m_n) \in \Gamma, 0 \leq m_i
< k_i\}$ a finite set.  To prove the lemma it is
sufficient to prove that each element of $\Gamma_0$ has a inverse in
$\Gamma$.  This is in view of (3.17).  Let $\underline{m}= (m_1,
\cdots, m_n) \in \Gamma_0$ be such that $m_1$ is non-zero and
minimal.  (if the all first co-ordinates if $\Gamma_0$ are zero we
take the next one).  We will first prove that there exist
$\underline{m}^1 \in \Gamma$ such that $m_1^1 = -m_1$.  Let $k_1=
\ell_1 m_1 + s_1, 0 \leq s_1 < m_1, \ell_1>0$.  Consider
$x=(\ell_1+1) \underline{m} - (k_1, \cdots,k_r) \in \Gamma$. Then
($\ell_1+1) m_1 - k_1 =- s_1 +m_1 \leq m_1$.  Further $0 < - s_1+m_1$.
Now we can assume that $x \in \Gamma_0$ without changing the first
co-ordinate in view of (3.17).  Now by minimality of $m_1$ it will
follow that $-s_1+m_1 = m_1$.  Now consider $\underline{m}+ (\ell-1)
\underline{m} - (k_1, \cdots, k_r) = (0, *,*, \cdots, )  \in \Gamma$.
 Take any $\underline{m}^1 \in \Gamma_0$.  Assume that $m_{1}^{1}
\neq 0$ write $m_1^1 = \ell_1^1 m_1 +s_1, 0 \leq s_1 < m_1$.  Arguing
as earlier we see that $s_1=0$.  Hence $y^1 = \underline{m}^1 -
\ell_1^1 \underline{m} = (0, *, \cdots,*) \in \Gamma$.  We can
further assume that $y^1 \in \Gamma_0$.  Thus we have  proved that
given $\underline{m} \in \Gamma_0$ there exists $\underline{m}^1 \in
\Gamma$ such that $\underline{m}+ \underline{m}^1 = (0, *, \cdots, *)
\in \Gamma_0$.  Repeating this process we get inverses of all
elements of $\Gamma_0$. [QED]. $\hfill{\Box}$

\paragraph*{(3.18) Theorem}   Let ${\cal G}$ be Kac-Moody Lie-algebra.
 Let $A= A_n=\C [t_1^{\pm}, \cdots, t_n^{\pm}]$ be  a Laurent
polynomial in $n$ commuting variables.  Let $\psi$ and
$\overline{\psi}$ as defined in (3.14) and (3.15) with each
$\lambda_i$ dominant integrable.  Then $V (\overline{\psi})$ is an
integral irreducible module for $\tilde{\cal G}_A$ with
finite dimensional weight spaces.  Further $V(\overline{\psi})$ is
isomorphic to the first component of $V(\psi) \otimes A$.

\paragraph*{Proof.}  From Lemma (3.16) $A_{\overline{\psi}}$ is
irreducible.  Hence $V(\overline{\psi})$ is a irreducible highest
weight module from Proposition (3.3).  From Proposition (3.5) it will
follow that $V(\overline{\psi})$ is isomorphic to the first component
of $V(\psi) \otimes A$.  Since each $\lambda_i$ is dominant integral
$V(\lambda_i)$ is integrable.  Since each $V(\lambda_i)$ is a highest
weight module it is known that $V= \displaystyle{\otimes_{i=1}^{N}}
V(\lambda_i)$ is module with finite dimensional weight space with
respect to $h$ where $h$ is included diagonally in $\oplus h$.  From
Lemma (3.11) it follows that the map defined in (3.10) is
surjective.  Hence $V \cong V (\psi)$.  Since $V(\psi)$ has finite
dimensional weight spaces it will follow from Lemma (3.6) that $V
(\overline{\psi})$ has finite dimensional weight spaces with respect
to $\tilde{h} = h \oplus D$. [QED].

\paragraph*{(3.19) ~ Remark}  Most often $A_{\overline{\psi}}=A$.  See
(1.7) of [E1].  In this case $V(\overline{\psi}) = V (\psi) \otimes
A$.

\paragraph*{(3.20) Proposition}  Suppose $V(\overline{\psi})$ is an
integrable highest weight module for $\tilde{\cal G}_A$ with finite
dimensional weight spaces.  Then
$$\overline{\psi} (h \otimes t^{\underline{m}}) =
\displaystyle{\sum_{j=0}^{N}} a_{I_j}^{\underline{m}} \lambda_j (h)
t^{\underline{m}} \ {\rm for} \ h \in h'$$
 for some distinct non-zero scalar
$(a_{i1, \cdots,} a_{iN_i})$ and $a_{I_j}^{\underline{m}}$ are as
defined in (3.10) and each $\lambda_j$ is dominant integral.

\paragraph*{Proof}  Let $\psi= E (1) \circ \psi$ where  $E (1)
t^{\underline{m}}=1$.  Let $v$ be a highest weight vector of $V
(\psi)$.  Then by Lemma (3.6), $V (\psi)$ has finite dimensional
weight spaces.  Then from the proof of Lemma (3.7),$ V (\psi)$ is a
module for ${\cal G}' \otimes A/I$ where $I$ is a  co-finite ideal of
$A$.  Further $I$ is generated by polynomials $P_j (t_j)$ for $1 \leq
j \leq n$.  We can assume that $P_j (t_j)$ has no zero roots as one
can multiply $P_j (t_j)$ by $t_j^{- \ell}, \ell >0$ and will get the
same ideal $I$.  Further we can assume that each polynomial $P_j$ is
not a constant.  In case $P_j$ is constant  then the module will be
trivial.  Let $a_{j1}, \cdots, a_{jN_j}$ be distinct non-zero
roots of $P_j (t_j)$.  Let
$$P_j^1 (t_j) = \displaystyle{\prod_{k=1}^{N_j}} (t_j- a_{jk}).$$

Let $I'$ be a co-finite ideal generated by $P_j^1 (t_j)$ insider $A$.

\paragraph*{Claim}  ${\cal G}' \otimes I^1$ is zero on $V(\psi)$.
Consider the map $\Phi: {\cal G}' \otimes A/I \to {\cal G}' \otimes
A/I' \to 0$.  It is easy to verify that ker $\Phi$ is solvable.  Let
$\alpha$ be a simple root of ${\cal G}$.  Let $X_\alpha, Y_\alpha,
H_\alpha$ be a $sl_2-$tripple and let ${\cal G}_\alpha$ be the span of
$X_\alpha, Y_\alpha$ and $H_\alpha$.  Consider $\overline{\cal
G}_\alpha = {\cal G}_\alpha \otimes A/I$ which is finite dimensional.
 Let $W$ be $\overline{\cal G}_\alpha $ module generated by $v$.
Since $V (\psi)$ is integrable and $X_\alpha \otimes P$ acts
trivially on $v$ and $H_\alpha \otimes P$ acts as scalars, we
conclude that $W$ is finite dimensional.  By restricting the action
of the solvable Lie-algebra ${\cal G}_\alpha \otimes I' / I$ to $W$
we get a vector $w$ in $W$ (by Lie-theorem) such that ${\cal
G}_\alpha \otimes I'/I$ acts as scalars on $w$. From the proof of
the Proposition (2.1) of [E3] (we assumed the module is irreducible but for this
conclusion we do not need it) we get that

\paragraph*{(3.21)} ${\cal G}_\alpha \otimes (I'/I) \  w=0$.

\paragraph*{Subclaim}  ${\cal G}_\alpha \otimes (I' /I) \ v=0$.  By
definition of $W$ we have a $Y \in U (Y_\alpha \otimes A/I)$ such that
$Y v =w$.  Since $V (\psi)$ is irreducible there exists $X \in U
(\tau)$ such that
$$X w = v = X Yv$$
where $X=X_-HX_+, X_+ \in U (N^+ \otimes A) X_- \in U (N^- \otimes
A),H \in U (h \otimes A)$.
First note that the weight of $V(\psi)$ are of the form $\lambda -
\sum n_i \alpha_i + \delta_{\underline{m}} $ where $\lambda= \psi
\mid h$ and $\alpha_1, \cdots, \alpha_n$ the simple positive roots
and $n_1, \cdots, n_a$ are non-negative integers.
$\delta_{\underline{m}}$ is a null root.  The weight of $w$ is
$\lambda - s \alpha$ for $s \geq 0$.  Thus for these reasons $X_+ \in
U (X_\alpha \otimes A)$.  Further $X_-$ has to be constant.  Thus $X$
is linear combinations of products of the form.

\paragraph*{(3.22)}  $$U (h \otimes A)X_{\alpha} Q_1 \cdots X_\alpha
Q_\ell, Q_i\in A.$$

First we will see that
\paragraph*{(3.23)} $$ {\cal G}_\alpha \otimes I' U (h \otimes A) w
=0 \ {\rm by} \ (3.21).$$
Now consider

\paragraph*{(3.24)}

\begin{eqnarray*}
{\cal G}_\alpha \otimes I' X_\alpha Q_1 w &=&
X_\alpha Q {\cal G}_\alpha \otimes I'w \\
&+& [X_\alpha, {\cal G}_\alpha] \otimes Q I' w =0
\end{eqnarray*}
Both terms are zero by (3.21).  Now by induction on $\ell$ we see
that
$${\cal G}_\alpha \otimes I' X_\alpha Q_1 \cdots X_\alpha Q_{\ell} w
=0.$$
>From (3.22), (3.23) and (3.24) it follows that ${\cal G}_\alpha
\otimes I' Xw = {\cal G}_\alpha \otimes I' v =0$.  This proves the
subclaim.  Since $h'$ is spanned by $H_\alpha, \alpha$ simple we have

$$h' \otimes I' v =0.$$

We will now show that $Y_\beta \otimes I'v=0$ for any positive root
$\beta$ where $Y_{\beta}$ is  a root vector of root $- \beta$.  We do
this by induction on the height of $\beta$.  We clearly know this for
$\beta$  such that height $\beta =1$ by (3.21).

Let $\alpha_1$ be any simple root.  Consider
$$X_\alpha Y_\beta \otimes I' v= Y_\beta \otimes I' X_\alpha
v+[X_\alpha, Y_\beta] \otimes I'v.$$
First term is zero since $v$ is a highest weight vector.  The second
term is zero by induction.  Thus we have proved that $Y_\beta \otimes
I'v$ is a highest weight vector in an irreducible highest weight
module $V(\psi)$.  Hence for weight reasons $Y_\beta \otimes I'v=0$.

Now consider
$$\tilde{W} =\{ w \in V (\psi), {\cal G}' \otimes I' w=0 \}$$
which is a non-zero submodule of $V(\psi)$.  Hence $\tilde{W}
=V(\psi)$.  This proves the claim.

Thus we have a module for ${\cal G}' \otimes A/I'$.  By by Lemma 3.11
(b) we have ${\cal G}' \otimes A/I' \cong {\cal G}'_N$ where $N=N_1
\cdots N_r$.  Consider $h'$ sitting in the ith place of ${\cal G}'_N$
which acts as linear function on the highest weight  vector $v$ of $V
(\psi)$.  By standard theory of integrable modules it follows that
$\lambda_i$ is dominant integral weight.  From the map given  in
Lemma 3.11 (a) it follows that
$$\psi (h' \otimes t^{\underline{m}}) = \sum a_{I_i}^{\underline{m}}
\lambda_i (h') \ {\rm for } \ h' \in h'.$$

\section*{Section 4}

We will now extract two classes of integrable irreducible modules for
toroidal Lie-algebra$\tau$ with finite dimensional weight spaces.

\paragraph*{(4.1) Example}  ${\cal Z}=0$ case.  Take ${\cal G}$ to be finite
dimensional simple Lie-algebra $\stackrel{\circ}{\cal G}$ in theorem
(3.18). Then $V (\overline{\psi})$ is an integrable irreducible
module with finite dimensional weight spaces for the toroidal
Lie-algebra $\tau$ where center ${\cal Z}$ acts trivially.

\paragraph*{(4.2) Example} ${\cal Z} \neq 0$ case.  Let $\stackrel{\circ}{\cal
G}$ be simple finite dimensional Lie-algebra.  Let ${\cal
G}_{af}=\stackrel{\circ}{\cal G} \otimes \C [t_1, t_1^{-1}] \oplus \C
C_1 \oplus \C d_1$ be the non-twisted affine Kac-Moody Lie-algebra.
Consider the
following Lie-algebra homomorphism
$$\Phi' : \stackrel{\circ}{\cal G} \otimes A \oplus {\cal Z}\oplus
\sum_{i=1}^{n} \C d_i \to {\cal G}'_{af} \otimes A_{n-1} \oplus
\sum_{i=1}^{n} \C d_i$$
by \\
(1) $ \Phi'$ is Id on $\stackrel{\circ}{\cal G}\otimes A_n$ \\
(2) $ \Phi' (d_i) = d_i$ \\
$$
(3)~~~ ~~~~~~~ \Phi' (t^{\underline{m}} k_i) =
\begin{cases}0  \ & {\rm if}   \ m_1 \neq
0 \\
0 & \ {\rm if}  \  m_1 =0,2 \leq i \leq n \\
C_1 \otimes t^{\underline{m}} \ & {\rm if} \ m_1=0, i=1
\end{cases}
$$
Here $A_{n-1} = \C [t^{\pm}_2, \cdots, t_n^{\pm}]$ and
$\stackrel{\circ}{\cal G} \otimes A$ is identified inside ${\cal
G}_{af} \otimes A_{n-1}$ by $X \otimes t^{\underline{m}} \to (X
\otimes t_1^{m_1}) \otimes t^{\overline{m}}$ where $\overline{m} =
(m_2, \cdots, m_n)$.  Now take ${\cal G}_{af}$ to be the
Kac-Moody Lie-algebra in Theorem (3.18) and take $n-1$ instead of
$n$. Then $V (\overline{\psi})$ is an irreducible integrable ${\cal
G}_{af} \otimes A_{n-1} \oplus \sum_{i=1}^{n} \C d_i$ module.
Then $V(\overline{\psi})$ is a module for $\tau$ via the homomorphism
$\Phi'$.

Towards the end of this paper we prove that the above two classes are
the only irreducible integrable modules with finite dimensional
weight spaces upto an automorphism of $\tau$.

We will now recall certain automorphisms of $\tau$ constructed in
section  4.3 of [E1].  Let $A = (a_{ij})_{\stackrel{1 \leq i \leq
n}{1 \leq j \leq n}}$ be an element of $GL(n, \Z)$ the group of
integral matrices of order $n$ with determinate $\pm 1$.  Let
$\underline{r} = (r_1, \cdots, r_n), \underline{s} = (s_1, \cdots,
s_n) \in \Z^n$.  Let $ e_i = (0, \cdots, 1, \cdots, 0)$ be such that
1 on the ith place and zero everywhere.  Let $A \underline{r}^T =
\underline{m}^T , A \underline{s}^T = \underline{d}^T$ where
$\underline{m}, \underline{d} \in \Z^n$.  Let $a (i) = (a_{1i},
\cdots, a_{ni})$ so that $A e_i^T =a (i)^T$.  $T$ denotes the transpose.  We
now define an automorphism of $\tau$ again denoted by $A$.

$$A(X \otimes t^{\underline{r}}) = X(\underline{m})  \ \ A(d
(t^{\underline{r}}) t^{\underline{s}}) = d (t^{\underline{m}})
\cdot t^{\underline{d}}$$

Let $(d_1^1, \cdots, d_n^1)=(A^T)^{-1} (d_1, \cdots, d_n)^T$.  Define $A
(d_i)= d_i^1$.  It is straightforward to check that $A$ defines an
automorphism of $\tau$.  A does not preserve the natural
$\Z^n-$gradation of $\tau$.

We note  that ${\cal Z}$  does not commute with $\tau$ but commutes with
$\stackrel{\circ}{\cal G} \otimes A \oplus {\cal Z}$.  In spite of this we
call them central operators since they are as good as central.

Let $V$ be irreducible $\tau$ module with finite dimensional weight
spaces.  We have the following Lemmas.

\paragraph*{ (4.3)  Lemma (1)}  Let $z \in
{\cal Z}$ be a homogeneous element of  degree $m$ such that $z v \neq 0$ for
some $v$ in $V$. Then $z w \neq 0$ for all non zero $w$ in $V$.  \\

(2) The zero degree central operators $K_i= C_i$ act by scalars

\paragraph*{Proof}  Let $W =\{v \in V; z v = 0 \}$.  Consider
$z X
\otimes t^{\underline{s}} v = X \otimes t^{\underline{s}} z v= 0$ for
$v$ in $W$.  Further $z d_i v = (d_i z+m_i z) v=0$ for $v$ in $V$.
Hence $W$ is a submodule of $V$.  Since $V$ is irreducible we should
 either have $V=W$ or $W=0$.  But we know that $W \neq V$ and hence
$W=0$.

(2)  The zero degree central operators $K_i$ leaves each finite
dimensional weight space invariant.  Thus $K_i$ has a eigenvector
$v$ of eigenvalue $\ell$.  Since $V$ is generated by $v, K_i$
(central) should act by the same scalar $\ell$ everywhere.  This
argument holds good for any operator on $V$ of zero degree which
commutes with $\tau$ and leaves a finite dimensional space invariant.
[QED].

We will now prove an important lemma which is crucial for the
classification result.

\paragraph*{(4.4.) Lemma}  Let $z$ be a central operator in $ {\cal Z}$ of
degree $\underline{m}$ such that $z \neq 0$.  Then there exists a
central operator $T$ (need not be in $Z$) on $V$ of degree
$-\underline{m}$ such that
$$Tz= zT=Id$$

\paragraph*{Proof} First consider $zV$ a non zero submodule of $V$.
Hence $zV=V$.  Now given a $v$ in $V$ there exist a unique $w$ in $V$
such that $z w=v$.  (uniqueness follow from the fact that $z w_1 =0$
implies $w_1=0$ by lemma 4.3 (1)).  Define
$T:V \to V$ such that $T(v) =w$. \\
Then clearly $z Tv =zw =v$.  Now consider
\begin{eqnarray*}
z (X (\underline{r}) Tv) &=& X (\underline{r}) zTv \\
&=& X (\underline{r}) v.
\end{eqnarray*}
Hence by definition $T(X (\underline{r})v) =X (\underline{r}) Tv$.
Thus $T$ commutes with ${\cal G}\otimes A$ and hence commutes with
$z$.  In particular
$$Tzv =zTv =zw =v. [QED]$$

\paragraph*{(4.5) Theorem}  Let $V$ be irreducible $\tau$ module with
finite dimensional weight spaces with respect to
$\stackrel{\circ}{h} \oplus Z_0 \oplus D$ where $ Z_0$ is the linear span of
$K_1, \cdots, K_n$.  Let
$$L=\{\underline{r},- \underline{r}; t^{\underline{r}} K_i \neq 0 \ {\rm on} \ V \ {\rm
for \ some} \ i\}.$$
Suppose the dimension of the $\C$ linear span of $L$ is $k$.  Then upto an
automorphism  of $\tau$

\begin{enumerate}
\item[(1)] there exists non-zero integers $m_1, m_2,
\cdots, m_k$ and nonzero central operators $z_1, \cdots z_k$ such
that the degree of $z_i$ is equal to $(0, \cdots, m_i, 0 \cdots
0)$
\item[(2)]  $k<n $
\item[(3)]  $t^{\underline{r}} K_i = 0$ for all $i$ and for all
$\underline{r}$ such that $r_j \neq 0$ for some $k+1 \leq j \leq n$.
\item[(4)] $t^{\underline{r}} K_i=0$ for $1 \leq i \leq k$ and for
all $\underline{r}$.
\item[(5)] There exists a proper submodule $W$ of
$\stackrel{\circ}{\cal G}_A \oplus {\cal Z}\oplus D_k$ such that $V/W$ has
finite dimensional weight spaces with respect to $h \oplus {\cal Z}_0+D_k$
where $D_k$ is the linear span of
\end{enumerate}
$$ \{d_{k+1}, \cdots, d_n\}.$$
\paragraph*{Proof}  Let $T = \{z_1 z_2 \cdots z_{\ell} \mid$
where $z_i =t^{ \underline{r}} K_j \neq 0$ or $z_i=$ inverse of
$t^{\underline{r}} k_j \neq 0$  for some $j$ and for some $\underline{r}\}$.

Let $L^1 = \{ \underline{s} \in \Z^n \mid $ where $\underline{s}$ is the
degree of some operator in $T \}$.  Then clearly $L^1$ is a sub group
of $\Z^n$.  Now by standard basis theorem  there exists $\underline{s}_1,
\cdots \underline{s}_n \in \Z^n$ such that  $m_1 \underline{s}_1
\cdots m_k \underline{s}_k$ is a basis of $L^1$ for some non-zero
integers $m_i$.  Now we can find an automorphism $B$ such that $B
\underline{s}_i = (0 \cdots 1, \cdots 0)$.  Let $\underline{r}_i =
m_i \underline{s}_i$ for $1 \leq i \leq k$.  Then $B \underline{r}_i =
(0, \cdots m_i, \cdots 0)$
for $i \leq i \leq k$.  So after twisting  the automorphism we can
assume that there exists non-zero central operators $z_1, \cdots z_k$,
of degree $(m_1, 0, \cdots 0) \cdots (0, \cdots m_k, 0 \cdots 0)$.
Let $T_1, \cdots T_k$ be the inverse
of $z_1, \cdots z_k$.

\paragraph*{Claim 1} $W_1= \{z_1 v-v \mid v \in V \}$ is a proper
$\stackrel{\circ}{\cal G}_A \oplus {\cal Z} \oplus D_1-$module.

Note that for $i \neq 1, d_i$ commutes with $z_1$ and hence $W_1$ is
$d_i-$invariant.  Clearly $W_1$ is $\stackrel{\circ}{\cal G}_A \oplus
{\cal Z}$
invariant.  To see that $W_1$ is proper, just note that $W_1$ cannot
contain $d_1$ weight vectors.  Now consider $W_2 = \{z_2 v-v; v \in
V/W\}$.  By above argument we can see that $W_2$ is a proper
$\stackrel{\circ}{\cal G}_A \oplus {\cal Z} \oplus D_2$ module. Continuing
this process we see that $W =\{z_i v-v \mid v \in V, 1 \leq i \leq k\}$
is a proper $\stackrel{\circ}{\cal G}_A \oplus {\cal Z}\oplus D_k-$module.

\paragraph*{Claim 2}  $T_i W=W$ and $v-T_iv \in W$ for all $v$ in
$V$.

Consider $v- T_i v = z_i T_i v -T_i v \in W$.  Now $T_i (z_j v-v) =
z_j (T_i v) - T_i v \in W$.  Hence the claim 2.

\paragraph*{Claim 3}  Let $\lambda $ be a weight of $V$.  Let
$\overline{V}_\lambda = \displaystyle{\oplus_{\underline{r}}}
V_{\lambda+ \delta_{\underline{r}}}$ where the sum runs over all
$\underline{r}$ such that $r_{k+1} =0 \cdots =r_n$.  Let
$\overline{W}_{\lambda} = \overline{V}_{\lambda} \cap W$.  Then $M=
\overline{V}_{\lambda} /_{\overline{W}_\lambda}$ is finite
dimensional.

To see this first recall that $z_i$ is central operator
of degree $\underline{m}_i = (0, \cdots, m_i, \cdots, 0)$.  Consider
\begin{eqnarray*}
S & =&\oplus_{\mid r_i \mid < \mid m_i \mid}
V_{\lambda+ \delta_{\underline{r}}} \\
r_{k+1} &=& 0 \cdots = r_n.
\end{eqnarray*}

Since each weight space is finite dimensional it follows that $S$ is
finite dimensional, so to prove the claim 3 it is sufficient to prove that
every $v$ in $V_{\overline{\lambda}}$ is in $S$ modulo
$W_{\overline{\lambda}}$.  Let $v \in
V_{\lambda+\delta_{\underline{r}}}, \underline{r} = (r_1, \cdots r_k,
0 \cdots 0)$.  Let
$$r_i = l_i m_i +g_i , \mid g_i \mid < \mid m_i)$$
for $g_i, l_i \in \Z$.  Note that the  sign of $l_i$ depends on the
sign of $r_i$ and $m_i$. Suppose $r_i > 0$ and $m_i <0$ then $l_i \leq
0$.  Suppose $r_i>0$ and $m_i<0$ then $l_i \geq 0$.

Let
$$
Z = \displaystyle{\prod_{l_i <0}} z_i^{-l_i}
\displaystyle{\prod_{l_i \geq 0}} T_i^{l_i}
$$

Then by claim 2 and the definition of $W$ we see that $Zv=v$ modulo
$W$ for all $v$ in $V$.  On the other hand it is easy to see that $Z
v \in S$. Hence the claim.

Thus we have proved that $V/W$ is a $\stackrel{\circ}{\cal G}_A
\oplus Z \oplus D_k-$module with finite dimensional weight space with
respect to $h \oplus Z_0 \oplus D_k$.  This provs (5).

(3)  Suppose $t^{\underline{r}} K_i \neq 0$ for some $\underline{r}$
and such that $r_j \neq 0$ for $k+1 \leq j \leq n$.  This increases
the dimension by 1.  So we have (3).

(1). Let $h_k$ be the linear span of $h \otimes t^{\underline{r}},
t^{\underline{r}} K_i$ for $1 \leq i \leq k$ and $r_j=0$ for $k+1
\leq j \leq n$ and $d_{k+1}, \cdots, d_n$.  One can check that $h_k$
is solvable.  In fact $[[h_k, h_k], h_k]=0$.  Now clearly
$V_{\overline{\lambda}}$ is a submodule for $h_k$ and
$W_{\overline{\lambda}}$ is a submodule.  Therefore $M$ is a finite
dimensional module for $h_k$.  Thus by Lie-theorem there exists a
vector $v$ in $M$ such that
$h (\underline{r}) v = \lambda (h, \underline{r}) v $ for $h
(\underline{r}) \in h_k$.  Consider for $\underline{m}, \underline{d}
\in \Z^n$ such that $d_j=m_j=0$ for $k+1 \leq j \leq n$.  Then
\begin{eqnarray*}
(h,h') d(t^{\underline{m}}) t^{\underline{d}} v &=& \left [h(\underline{m}), h'
(\underline{d})\right ] v \\
&=& (h (\underline{m}) h' (\underline{d})-h (\underline{m}) h
(\underline{d})) v \\
&=&(\lambda (h, \underline{m}) \lambda (h', \underline{d}) - \lambda
(h',\underline{d}) \lambda (h, \underline{m})) \\
&=& 0
\end{eqnarray*}

Let $z=d (t^{\underline{m}}) t^{\underline{d}}$ and suppose $z \neq
0$ on $V$.  Then we have proved that there exists a non-zero vector
in $V/W$ (there is $v \notin W)$ such that $z v \in W$.  But
$v=z^{-1} z v \subseteq z^{-1} W \subseteq W$.  A contradiction.
Hence $z=0$.  In particular $t^{\underline{r}} K_i =0$ for $1 \leq i
\leq k$ and $r_j =0$ for $k+1 \leq j \leq n$.  (see 2.2 (1)).  This
together with (3) proves (4).  The second part of (1) follows from
above.

(5) Suppose $k=n$.  Then by 4, ${\cal Z}$ has to be zero a contradiction to
the fact that $k=n>0$.  [QED].

We record here the following Lemma about the
dimensions of central operators acting on $V$.  We will not need it
anywhere but of independent interest.

\paragraph*{(4.6) Lemma}  Let $z_1, z_2$ be non-zero central
operators of degree $\underline{m}$ of an irreducible module $V$ of
$\tau$ with finite dimensional weight space.  Then $z_1=kz_2$ for
some constant $k$.  In particular there is at most one dimensional
non-zero central operator in a given degree.

\paragraph*{Proof}  Let $T_1$ be the inverse central operator
 of degree $- \underline{m}$ of $z_1$.  Consider $z_2 T_1$ which is of
degree zero central operator and leaves a finite dimensional space
$V_{\underline{g}}$ (see the proof of Lemma 2.6) invariant.

Thus it has a eigenvector $v$ say of eigenvalue $k$.  But $v$
generates $V$ and hence $z_2 T_1= k$ on $V$.  Now $z_1 k = z_1 z_2
T_1=z_2$.  Hence we are done. [QED].

We need one more reduction modulo an automorphism of $\tau$ before we
can take up the classification
problem.  Recall that $A \in GL(n, \Z)$ defines an automorphism of
$\tau$ such that $A d(t^{\underline{r}}) t^{\underline{s}} = d (t^{A
\underline{r}}). t^{A \underline{s}}$.  (We  are supressing Transpose
$T$ and there is no confusion).  It is easy to see that if $A =
(a_{ij})$ then

\paragraph*{(4.7)} $A(t^{\underline{s}} K_i) =
\displaystyle{\sum_{j=1}^{n}}  a_{ij} t^{A (\underline{s})} K_j$.  By
taking $\underline{s}=0$ we have $A (K_i) =
\displaystyle{\sum_{j=1}^{n}} a_{ij} K_j$.  Let $K_i$ act on $V$ by
$k_i$.  We know that there exists $k<n$ such that $k_i=0$ for $1 \leq
i \leq k$ upto an automorphism of $\tau$.  Now choose $A \in GL(n,
\Z),  A= \begin{pmatrix}I & 0 \\ 0 & B \end{pmatrix}$ such that $I$
is identity matrix of order $k \times k$ and $B \in GL(n-k, \Z)$ such
that $B (k_{k+1}, \cdots, k_n)^T = (0, \cdots, 0, \ell)$.

\paragraph*{(4.8) Proposition}  Let $V$ be irreducible module for
$\tau$  with finite dimensional weight spaces.  Let $k$ be an integer as
defined in Theorem (4.5).  Then upto automorphism of $\tau$  the
assertions of Theorem 4.5 holds and further we can assume $K_i =0$
for $ 1 \leq i \leq n-1$.

\paragraph*{Proof}  We have Theorem (4.5).  Now choose $A$ such that
(4.7) holds.

\begin{enumerate}
\item[(1)] Clearly holds.
\end{enumerate}
\paragraph*{(3)~ Claim} $t^{A (\underline{r})} K_i=0$ for all $i$ and $j$
such that $A(\underline{r})_j \neq 0$ for some $j$ such that
$k+1 \leq j \leq n$.  First
note that $r_{i} \neq 0$ for some $k+ 1 \leq i \leq n$ by the choice
of $A$.  Now by Theorem 4.5 (3) we have $t^{\underline{r}} K_i=0$ for
all $i$.  Then from (4.7) and the fact that $A$ in invertible the
claim follows.

\paragraph*{(4)~ Claim} $t^{A(\underline{r})} K_i =0 \ 1 \leq i
\leq k$ for all $\underline{r}$.  From (4.7) we have the claim by
noteing that $a_{ij} =0$ for $j \geq k+1$.

2 and 5 follows from arguments similar to Theorem 4.5 (2) and (3).
In addition we can assume that $A (K_i)=0, 1\le i\le n-1$.  [QED].

We will now start eleminating several cases in order to classify
integrable irreducible modules for $\tau$ with finite dimensional
weight spaces.

\paragraph*{(4.9) Proposition}  Let $V$ be irreducible integrable
module for $\tau$ with finite dimensional weight spaces.   Suppose $k
< n-1$ (see Theorem 4.5) and suppose $k_i \neq 0$ for some $i$,  Then
such $V$ does not exists.

\paragraph*{Proof}  By Proposition (4.8) we can assume that $k_i=0$
for $1 \leq i \leq n-1$ and $k_n \neq 0$.  Let ${\cal G}_{af}=
\stackrel{\circ}{{\cal G}} \otimes \C [t_n, t_n^{-1}] \oplus \C K_n
\oplus \C d_n$ be an affine Lie-algebra.  Let ${\cal G}_{af} =
N^+ \oplus h' \oplus N^-$ be a standard decomposition with $h'=
\stackrel{\circ}{h} \oplus \C K_n \oplus \C d_n$.  Now by Proposition
(2.4) we can assume that there exists a weight vector $v$ of weight
in $V$ such that
$$N^+ \otimes A_{n-1} v=0.$$
(The other case can be dealt similarly). Here $A_{n-1} = \C
[t_1^{\pm}, \cdots, t_{n-1}^{\pm}]$.  (We have chosen $n$ instead
of 1 in Proposition (2.4)).  In particular $h \in
\stackrel{\circ}{h}$ we have $h \otimes t^{\ell}_{k+1} t^m_n v =
0$ for all $m >0$ and for all $\ell$.

\paragraph*{Claim}  For any $0 \neq h \in \stackrel{\circ}{h}$ the
following vectors are linearly independent in $V$,  Fix $m >0$ and
$\ell$.
$$\{ht^{\ell - d}_{k+1} t^{-m}_{n} \cdot h t^d_{k+1} t_n^{-(m+1)}
v, d \in \Z \}.$$
Suppose there exists non-zero scalars $a_d$ such that
$$B = \sum a_d ht^{\ell-d}_{k+1} t^{-m}_{n} ht^d_{k+1} t^{-
(m+1)}_{n} v=0.
$$
  Choose $h' \in h$ such that $(h, h') \neq 0$.
Consider \\
$h' t_{k+1}^s t_n^{m+1} B =0$ for any $s$, which implies
\begin{eqnarray*}
&&\sum a_dh t_{k+1}^{\ell-d}t_n^{-m} h' t^s_{k+1} t_n^{m+1}ht^d_{k+1}
t_n^{-(m+1)}v \\
&+& (h,h') \sum a_d s t_{k+1}^{\ell-d+s}t_n K_{k+1} h t^d_{k+1}
t_n^{-(m+1)}v \\
&+& (h,h') \sum a_d (m+1) t_{k+1}^{\ell-d +s} t_n K_n ht^d_{k+1}
t_n^{-(m+1)}v
\end{eqnarray*}
(by 2.2 (1)).

The term $t_{k+1}^{\ell-d+s} t_n K_{k+1}$ and $t_{k+1}^{\ell -d+s}
t_n K_n$ are zero (by Theorem 4.5 (3)) being central. Thus the second
and third term above are zero.

The first term is equal to
\begin{eqnarray*}
0&=&\sum a_dh t_{k+1}^{\ell-d}t_n^{-m} h t^d_{k+1}
t_n^{-(m+1)}h't^s_{k+1} t_n^{m+1} v \\
&+&(h,h') \sum a_dh t_{k+1}^{\ell-d} t_n^{-m} s t^{s+d}_{k+1} K_{k+1} v
\\
&+& (h,h') \sum a_dh t_{k+1}^{\ell-d}t_n^{-m} (m+1) t^{s+d}_k K_n v
\end{eqnarray*}

(by 2.2) (1)).

Fix a $d_0$ in the above and let $s=-d_0$.  Now from  (2.1)
$t^{s+d}_{k+1} K_{k+1} =0$ for $s+d \neq 0$.  And $K_{k+1}=0$ by
Proposition (4.8).  Thus the second term is zero in the above.  Now
by Theorem 4.5 (3), $t^{s+d}_{k+1} K_n=0$ for $s+d \neq 0$.  Hence
the third term is zero but for $d_0$.  The first term is zero being
an highest weight.  Hence we have
$$a_{-s} ht^{\ell+s}_{k+1} t^{-m}_n \cdot K_n v =0$$
\paragraph*{(4.10)}
Suppose $ht^{\ell+s}_{k+1} t^{-m}_{n} v =0$.
Consider
\begin{eqnarray*}
0&=& h' t^{- (\ell+s)}_{k+1} t^m_n ht^{\ell+s}_{k+1} t^{-m}_{n} v
\\
&=& - (h, h') (\ell+s) K_{k+1} v \\
&+& m (h, h') K_n v
\end{eqnarray*}
(by  2.2 (1))  and $v$ is a highest weight.  We know that $K_{k+1} v
=0$.  But the second term is non zero.  Thus (4.10) is false.   This
proves $a_{-s}=0$ a contradiction.  This proves our claim.  Hence we
have proved that under the conditions of Proposition (4.9).
$V_{\lambda+ \ell \delta_{k+1} - (2m+1)\delta_n}$ is infinite
dimensional.  This completes the proposition. [QED].

\paragraph*{(4.11) Remark} The case $k=n-1$ where modules exists (see
Examples 4.2) will be dealt in the next section.  The case $k=0$ and
$K_i=0$ for all $i$, in which case ${\cal Z}=0$ is dealt in [E3].

Now we will deal the remaining case where $k \geq 1$ and $K_i=0$ for
all $i$.

We will first recall an important result due to Futorny [F] on
Hiesenberg Lie- Algebra.  Let $H$ be a finite dimensional vector
space with non-degenerate symmetric billinear form (,).  Then $L(H) =
H \otimes \C [t, t^{-1}] \oplus \C c $ is called Hisenberg
Lie-algebra with the following bracket
$$[h \otimes t^{m}, h' \otimes t^{\ell} ] = (h,h') m \delta_{m+ \ell,
0} C.$$

\paragraph*{(4.12) Proposition}  (Proposition 4.3 (i) [F]).  Let $V$
be any $\Z-$graded.  $L(H)-$module with finite dimensional graded
spaces and center acts by nonzero scalars. Then $V$ admits a
graded vector $v$ such that $H \otimes t^n v=0$ for all $n>0$ (or
for all $n<0)$.

\paragraph*{Proof}  It is only proved for one dimensional $H$.  But
the proof works for any finite dimensional $H$ by choosing orthogonal
basis for $H$.  There it is assumed that $V$ is irreducible.  But it
is not needed  for geting a highest weight (or lowest weight) vector.
[QED].

\paragraph*{(4.13) Proposition.}  Let $V$ be integrable irreducible
module for $\tau$ with finite dimensional weight spaces.  Let $k$ be
as in theorem 4.5.  Suppose $k \geq 1$ and $K_i=0$ for all $i$.  Then
such a module $V$ does not exists.

\paragraph*{Proof} Recall that $\stackrel{\circ}{\cal G}$ is finite
dimensional simple Lie algebra with $\stackrel{\circ}{\cal G}=n^+
\oplus \stackrel{\circ}{h} \oplus n^-$.  Then by Proposition (2.12)
there exists a weight vector $v$ of $V$ such that $n^+ \otimes A_n
v=0$.  (The case $n^- \otimes A_n v=0$ can be done similarly).
Suppose $t^{\underline{m}} K_i \neq 0$.  Then by Theorem 4.5 (3), (4)
we have $i \geq k+1$ and $m_{k+1} =0 \cdots = m_n$.  Let $H$ be Hisenberg
Lie-algebra spanned by $ht^{\underline{m}} t^k_n, k >0, ht^{-k}_n,
k>0$ and $t^m K_i$ with Lie-braket.
$$[h t^{\underline{m}} t^k_i, h' t^{+ \ell}_k] = (h, h')
kt^{\underline{m}} K_i \delta_{k+\ell
0}.$$
  Consider $M$ the $H$ module generated by $v$.  Then by
Proposition (4.12) there  exists  $w$ in  $M$
such that
$$ (a) ~~~~~~ ht^m t_i^k  w =0, k>0 \ {\rm or}  \ (b) ~~~~ ht^ k_i w =0 \
{\rm for } \ k <0.$$
Assume (a) \\
Now $w =X v \ {\rm for } \ X \in U (h_A)$.  Then it is easy to see
that $n^+ \otimes A_n w=0$.  Let $\lambda$ be the weight of $w$.

\paragraph{Claim}  $\lambda \mid \stackrel{\circ}{h} \neq 0$.
Suppose it is zero.  Let $\alpha$  be a simple root in
$\stackrel{\circ}{\triangle}$ and let $X_{\alpha}, Y_{\alpha}, h_{\alpha}
=[X_\alpha, Y_\alpha]$ be as $sl_2-$ copy  inside
$\stackrel{\circ}{\cal G}$.  Then $X_\alpha \otimes
t^{\underline{s}}, Y_\alpha \otimes t^{-s}, h_\alpha$ is an $sl_2$
copy.  (because there is no zero degree centre).  Now by $sl_2$
theory for integrable module $Y_{\alpha} \otimes t^{-s} w =0$ for
any $\underline{s}$ as $\lambda (h_\alpha)=0$.  Thus we have
$\stackrel{\circ}{\cal G}\otimes A_n w =0$.  This implies ${\cal Z}$ is
zero which is not the case.  Thus $\lambda \mid \stackrel{\circ}{h}
\neq 0$.  Let $\alpha$ be a simple root such that $\lambda (h_\alpha)
 \neq 0$.  Let $X_\alpha =X, Y_\alpha =Y \ h=h_\alpha$ be an
$sl_2$ copy.  That is $[X,Y]=h, [h,Y]=-2Y, [h,X]=2X$.  Let $\ell_0$ be
such that

\paragraph*{(4.14)} $ht^{\underline{m}} - \ell_0 (X,Y)
t^{\underline{m}} K_i=0$.
where  $X t^{\underline{s}}$ we mean $X \otimes t^{\underline{s}}$.

\paragraph*{Claim}  The following infinite set of vectors are
linearly independent in $V$.
$$\{Yt^{-r}_{i} Yt^r_i w, r >0, r \neq \ell_0, -\ell_0\}.$$
Suppose there exists non-zero scalars $a_r$ such that
\paragraph*{(4.15)} $\sum a_r Y t^{-r}_{i} Y t^r_i w =0$.  We will be
using the following in the calculation below.
\begin{enumerate}
\item[(1)] $w$ is a highest weight vector
\item[(2)] $t_i^{\ell} K_i=0$ which follows from definition 2.2. for
$\ell \neq 0$ and by assumption for $\ell=0$.
\item[(3)] $d (t^{\underline{m}} t_i^{\ell}) t_i^{-r} =r
t^{\underline{m}} t_i^{\ell-r} K_i $ (by (2.2) and (2.1) and
$m_i=0)$.
\end{enumerate}
Consider for $r > 0$

\begin{eqnarray*}
&&X t_i^s X t^{\underline{m}} t_i^{\ell} Y t_i^{-r} Y t_i^{r} w \\
&=& ht^{s-r}_{i} h t^{\underline{m}} t_i^{\ell+r} w \\
&&-r (X,Y)h t_i^{s-r} t^{\underline{m}} \cdot t_{i}^{\ell+r} K_i w
\\
&+& h t_i^{r+s} h t^{\underline{m}}_{\ell+s} t_i^{\ell-r} w \\
&& -2 ht^{\underline{m}} t_i^{\ell+s}  w \\
&&-2 s (X,Y) t^{\underline{m}} t_i^{\ell+s} w \\
&+& r (X,Y) t^{\underline{m}} t_i^{\ell-r} K_i h t_i^{r+s} w
\\
&-& r (X,Y)^2 s t_i^{r+s} K_i t^{\underline{m}} t_i^{r+l} K_iw
\\
&+& ht_i^{r+s} h t^{\underline{m}} t_i^{-r+l} w
\end{eqnarray*}

Choose $\ell$ such that $\ell-r>0$ which implies $\ell+r >0$. Then
first three terms, the sixth, seventh and eighth  term are zero.
Thus applying $X t_i^{s}X t^m t_i^{\ell}$ to (4.15) such  that
$\ell -r
>0$. We have for all $r$
$$\sum a_r (ht^{\underline{m}} t_i^{\ell+s} + s (X,Y) t^m
t_i^{\ell+s}) w=0$$
choose $s$ such that $\ell+s=0$.  From 4.14 it follos that
\paragraph*{(4.16)}

$\sum a_r =0$.
Now choose $r_0$ be the maximal among $r$ that occur in (4.15).  Now
choose $\ell$ such that $\ell-r_0=0, \ell-r >0$ for $r \neq r_0$.
This implies $\ell+r>0 \ \forall  r$.  Again apply $X t_i^s
Xt^{\underline{m}} t_i^{\ell}$ to (4.15) we have

$$a_{r_0} ht_i^{r_0+s} ht^{\underline{m}} w -2 \sum a_r
(ht^{\underline{m}} t_i^{\ell+s}+s (X,Y)
t^{\underline{m}} t_i^{\ell+s} K_i)w$$
$$+r_0 a_{r_0} (X,Y) t^{\underline{m}} K_i ht_i^{r_0+s}w. $$
Now choose $s$ such that $r_0+s=0$.    Then $a_{r_0}
(ht^{\underline{m}}+ r_0 (X,Y) t^{\underline{m}} k_i) h w=0$.  By
choice of $r$ and the fact that $hw= \lambda (h) w \neq 0$, we
conclude that $a_{r_0}=0$.  A contradiction.  Thus $V_{\lambda+2
\alpha}$ is infinite dimension.

(b) $ht_i^k w=0, \ \ k < 0$.  Consider the set $Y t_i^r Yt_i^{-r}, r>0$
and apply $Xt^{\underline{m}} t_i^s X t_i^{\ell}$ then we get the
desired linearly independent set. [QED].

\section*{5. Section}  In this section we will deal with the last
case $k=n-1, K_n \neq 0$ and $K_i =0$ for $1 \leq i \leq n-1$.  We
will prove that in this case we will get all the modules defined in
example (4.2).

\paragraph*{(5.1) Lemma}   Any $\Z^{n-1}$ graded simple, commutative
Algebra $M$ such that each graded component is finite dimensional
over $\C$ is isomorphic to a subalgebra of $A_{n-1}$ such that each
homogeneous element is invertible.

\paragraph*{Proof}  Let $M = \displaystyle{\oplus_{\underline{r} \in
\Z^{n-1}}} M_{\underline{r}}$ where each $M_{\underline{r}}$ is
finite dimensional.  Since $M$ is graded simple it follows that $M_0$
is simple commutative Algebra of finite dimension over $\C$.  Then
clearly $M_0 \cong \C$.  Let $0 \neq w \in M_{\underline{r}}$.  The
ideal generated by $w$ has to be $M$ and hence there exists inverse
say $w^{-1}$.  Consider $w^{-1} M_{\underline{r}} \subseteq \C$.
Since $w^{-1} w =1$ it follows that $w^{-1} M_{\underline{r}} =\C$.
In particular each non-zero $M_{\underline{r}}$ is $\C$.  The Lemma
follows. [QED].

\paragraph*{(5.2) Theorem}  Let $V$ be irreducible integrable models
for $\tau$ with finite dimensional weight spaces.  Let $k$ be as
defined in Theorem (4.5).  Assume $k=n-1$ and $K_n \neq 0$.  Then $V$
is isomorphic to $V(\overline{\psi})$ as defined in Example (4.2).

\paragraph*{Proof}  First note that $t^{\underline{r}} K_i =0$ for
all $\underline{r}$ and $1 \leq i \leq n-1$ and $t^{\underline{r}} K_n
=0$ for all $\underline{r}$ such that $r_n \neq 0$.  Write

\paragraph*{(5.3)} $\tau = N^- \otimes A_{n-1} \oplus
\stackrel{\circ}{h} \otimes A_{n-1} \oplus \displaystyle{\sum_{r_n
=0}} t^{\underline{r}} K_n \oplus D + {\cal Z}' \oplus N^+ \otimes A_{n-1}$
where ${\cal Z}'$ is
spanned by
$t^{\underline{r}} K_i$ for $1 \leq i \leq n-1$ for all
$\underline{r}$ and $t^{\underline{r}} K_i$ for all $\underline{r}$
such that $r_n \neq 0$.  We have $N^- \oplus \stackrel{\circ}{h} \oplus
\C K_n \oplus \C d_n \oplus N^+={\cal G}_{af} =
\stackrel{\circ}{\cal G}\otimes \C[t_n, t_n^{-1}] \oplus \C K_n
\oplus \C d_n$.  We have already noted that ${\cal Z}'$ acts trivially on
$V$.  By Lemma 2.3 (5) we know that $K_n$ acts by an integer.  By an
automorphism we can assume that $K_n$ acts by a positive integer.
Now by Proposition (2.4) we get a weight vector $v$ in $V$ such that

\paragraph*{(5.4)}  $N^- \otimes A_{n-1} v=0$.  Let $\underline{h''}$
be the abelian Lie-algebra spanned by $\stackrel{\circ}{h} \otimes
A_{n-1} \oplus \displaystyle{\sum_{r_n=0}} t^{\underline{r}} K_n$.
Here $A_{n-1} = \C [t_1^{\pm}, \cdots, t_{n-1}^{\pm}]$.  Let $M$ be a
$\underline{\underline{h}}^{''} \oplus D-$ module generaed by $v$.

\paragraph*{Claim}  $M$ is irreducible
$\underline{\underline{h}}^{''} \oplus
D-$module.  Let $w$ be a weight vector of $M$.  Since $V$ is
irreducible there exists $X \in U (\stackrel{\circ}{\cal G}_A \oplus
{\cal Z}')$ such that $X w =v$.  By PBW theorem
$$X=X_- HX_+ \ {\rm where} $$
$X_+ \in U (N^+ \otimes A_{n-1}), X_- \in U (N^- \otimes A_{n-1})$
and $H \in U (h'')$ from (5.3).  $D$ and ${\cal Z}'$ does not appear
as $w$ is a
weight vector and ${\cal Z}'$ acts trivially.  But $Y v=w$ for some $Y \in U
(\underline{\underline{h}}^{''})$ we can see that $N^+ \otimes A_{n-1} Yv =0$.
Thus $X_+ w=0$, which means $X_+$ cannot appear and by weight reasons
$X_-$ cannot appear.  Thus $X=H$ which belongs to
$U(\underline{\underline{h}}^{''})$.  This proves the claim.

Since $v$ is a weight vector, $D$ acts by scalar and hence
$M=U(\underline{\underline{h}}^{''})v$.  So we have $M $ a $\Z^{n-1}$ graded
irreducible module for $ \underline{\underline{h}}^{''}$.  Thus $M \cong U
(\underline{\underline{h}}^{''}) /I$ for some graded Ideal $I$.  Since
$\underline{\underline{h}}^{''}$ is abelian $M$ is $\Z^{n-1}$ graded simple
commutative algebra.  Now by Lemma (5.1), $M$ is isomorphic to a
subalgebra of $A_{n-1}$.  Let $\overline{\psi}: U
(\underline{\underline{h}}^{''}) \to M= A_{\overline{\psi}} \subseteq
A_{n-1}$ be the quotient map
which $\Z^{n-1}$ graded.  By Proposition 3.3, we have $V \cong V
(\overline{\psi})$.  As earlier $\psi= E (1) \circ \overline{\psi}$
where E(1) $t^{\underline{m}} =1$.

$h'= \stackrel{\circ}{h} \oplus \C K_n$ and identify
$t^{\underline{m}} K_n $ as $K_n \otimes t^{\underline{m}} (m_n=0)$.
Thus $h' \otimes A_{n-1} = \stackrel{\circ}{h} \otimes A_{n-1} \oplus
\displaystyle{\sum_{m_n=0}} \C t^{\underline{m}} K_n$.  From
Proposition (3.20) we have
$$\psi (h \otimes t^{\underline{m}})=\displaystyle{\sum_{j=0}^{N}}
a_{I_J}^{\underline{m} } \lambda_j (h) \ {\rm for} \ h \in h'$$
and $\underline{m} \in \Z^{n-1}$.  Now extend each $\lambda_j$ to
$h$ that is give some  value to $\lambda_j (d_n)$ such that $\psi
(d_n) = \sum \lambda_j (d_n)$.  It is well known that the irreducible
integrable highest weight module $V(\lambda_j)$ for ${\cal G}_{af}$
are all isomorphic for various values of $\lambda_j (d_n)$.  (see
[K]).  Now $\displaystyle{\otimes_{j=1}^{N}} V (\lambda_j) \cong
V(\psi)$ being the unique irreducible module for the highest weight
$\psi$. [QED].

\paragraph*{Remark ~~ (5.5)}  We need to consider more general
modules than in [E3] where center acts trivially like in our  Example
(4.1).  Nevertheless we proved that an irreducible integrable module
(not graded) $\stackrel{\circ}{\cal G}_A$ is actually a module for
$\stackrel{\circ}{\cal G}\otimes A /I ($ Lemma (1.2)
and Proposition (2.1) of [E3]) where $I$ is a cofinite ideal of A
generated by polynomials with distinct roots.  But our lemma 3.11 (b)
says that $\stackrel{\circ}{\cal G}\otimes A/I \cong \oplus
\stackrel{\circ}{\cal G}$.  Thus in this case all irreducible
integrable modules with finite dimensional weight spaces for
$\stackrel{\circ}{\cal G}_A$ are given in Example (4.1).

\paragraph{Acknowledgements}  I record my sincere thanks to Fields
Institute's (Toronto) hospitality during the fall of 2000 where some of
the work has been done.  I also thank V. Futorny for bringing into my
notice the reference [F].

\newpage

\begin{center}
{\bf References}
\end{center}
\begin{enumerate}
\item[{[BB]}] Berman, S. and Billig, Y. {\it Irreducible representations for
toroidal Lie-algebras}. Journal of Algebra, 221, 188-231 (1999).
\item[{[BC]}] Berman, S and Cox, B. {\it Enveloping Algebras and
Representations of toroidal Lie-algebras}.  Pacific Journal of
Mathematics, 165 (2), 239-267 (1994).
\item[{[BGK]}] Berman, S., Gao, Y., Krylyuk, Y.: {\it Quantum tori
and the structure of elliptic quasi-simple Lie-Algebas,} Journal of
Functional Analysis, 135, 339-389 (1996).
\item[{[C]}] Chari, V. {\it Integrable representations of Affine
Lie-algebras.}  Invent Math. 85, 317-335 (1986).
\item[{[CP]}] Chari, V. and Pressley, A.N. {\it New Unitary
Representations of Loop Groups}.  Math. Ann. 275, 87-104 (1986).
\item[{[E1]}] Eswara Rao, S. {\it Iterated Loop Modules and a
filteration for Vertex Representation of toroidal Lie-algebras}.
Pacific Journal of Mathematics, 171 (2), 511-528. (1995).
\item[{[E2]}] Eswara Rao, S. {\it Classification of Loop Modules with
finite dimensional weight spaces}. Math. Ann. 305, 651-663 (1996).
\item[{[E3]}]  Eswara Rao, S. {\it Classification of irreducible
integrable modules for multi-loop algebras with finite dimensional
weight spaces.}   Journal of Algebra, 246, 215-225 (2001).
\item[{[E4]}]  Eswara Rao, S. {\it A generalization of irreducible
modules for toroidal Lie-algebras} TIFR preprint (2000).
\item[{[EM]}] Eswara Rao, S and Moody, R.V. {\it Vertex
representations for $n-$toroidal Lie-algebras and a Generalization of
the Virosoro Algebra}, Communications of Mathematical Physics, 159,
239-264 (1994).

\item[{[F]}] Futorny, V. {\it Representations of affine Lie-algebras,}
Queen's papers in Pure and Applied Mathematics,  (1997).
\item[{[H]}] Humphreys, J.E. {\it Introduction to Lie-algebras and
representations theory.}  Springer, Berlin, Hidelberg, New York (1972).
\item[{[K]}] Kac, V. {\it Infinite dimensional Lie-algebras}, Cambridge
University Press, Third Edition (1990).
\item[{[MEY]}] Moody, R.V., Eswara Rao, S. Yokomuma, T. {\it Toroidal
Lie-algebra and Vertex Representations}, Geom, Ded., 35, 283-307 (1990).
\item[{[MS]}] Moody, R.V., Shi, Z. {\it Toroidal Weyl groups}, Nova
Journal of Algebra and Geometry 1, 317-337 (1992).
\item[{[YY]}]  Youngsun Yoon, {\it On the polynomial representations of
Current Algebras}, Yale Prepprint (2001).
\end{enumerate}
\vskip 5mm
School of Mathematics \\
Tata Institute of Fundamental Research \\
Homi Bhabha Road \\
Mumbai 400 005 \\
India \\ [5mm]
e-mail: senapati@math.tifr.res.in

\end{document}